\documentclass{article}
\usepackage{amssymb}
\usepackage{amsmath}
\usepackage{graphicx}
\newtheorem{theorem}{Theorem}
\newtheorem{definition}{Definition}
\newtheorem{lemma}[theorem]{Lemma}
\newtheorem{rem}[theorem]{Remark}

\newcommand{\dom }{\,{\rm dom}\,}

\newcommand{\cover}[1]{\stackrel{#1}{\Longrightarrow}}
\newcommand{\comment}[1]{\mbox{}}

\def\qed{{\hfill{\vrule height5pt width3pt depth0pt}\medskip}}

\begin{document}
\begin{center}
   { \LARGE
Computer assisted proof of the existence of homoclinic tangency for the H\'enon map and for the
  forced-damped pendulum } \\
\vskip 0.5cm
  \vskip\baselineskip
    {\large
    Daniel Wilczak
    } \\
    Departament of Mathematics, University of Bergen, \\
Johannes Brunsgate 12, 5008 Bergen, Norway\\
and \\
Jagiellonian University, Institute of Computer Science, \\
 {\L}ojasiewicza 6, 30-348 Krak\'ow, Poland \\
\vskip 0.5cm
    {\large
    Piotr Zgliczy\'nski\footnote{
            Research supported in part by Polish State Ministry of
            Science and Information Technology  grant N201 024 31/2163 
         }
    } \\
 Jagiellonian University, Institute of Computer Science, \\
 {\L}ojasiewicza 6, 30-348 Krak\'ow, Poland
\vskip 0.5cm
  e-mail: wilczak@ii.uj.edu.pl, zgliczyn@ii.uj.edu.pl
\\
\vskip 0.5cm
    \today
\end{center}

\begin{abstract}
We present a topological method for the efficient computer
assisted verification of the existence of the homoclinic tangency
which unfolds generically in a one-parameter family of planar
maps. The method has been applied to the H\'enon map and the
forced damped pendulum ODE.
\end{abstract}


\section{Introduction}
The goal of our paper is to describe a  method for the
verification of the existence of the quadratic homoclinic tangency
which unfolds generically in a one-parameter family of planar
maps. This is an important problem problem in dynamics, because
establishing the existence of the generic homoclinic tangency has
very deep dynamical consequences, see for example
\cite{PT,WY1,WY2,RC} and  references given there.

Our paper was inspired by the works of Arai and Mischaikow
\cite{A,AM}, who combined some tools from  the Conley index
theory, the tools of computational homology  from \emph{the CHomP
project} \cite{CHOMP} and the set oriented numerical methods from
\emph{the GAIO project} \cite{GAIO}, into a method for computer
assisted proof  of the existence of the generic homoclinic
tangency. Using it they proved in \cite{AM} the existence of the
generic homoclinic tangency for the dissipative H\'enon map
$H_{a,b}(x,y) = (a-x^2 +by,x)$ for parameter values close to $a =
1.4$, $b = 0.3$ and $a = 1.3, b = -0.3$. Their method contains
essentially two separate parts: first using the Conley index
approach they prove the existence of the homoclinic tangency for
some parameter value (steps 1 to 5 in the terminology used in
\cite{AM}) and in the second part (step 6 in \cite{AM}) they
verify some transversality-type condition, which implies the
genericity of the homoclinic tangency established in the first
part. The computation times reported in \cite{AM} are around 260
minutes and 100 minutes on a PowerMac G5 (2GHz) for $b=0.3$ and
$b=-0.3$, respectively. In these computations the second part took
61 minutes and 24 minutes, respectively. From these computation
times it is quite clear that there is a little hope to apply
successfully this method to ODEs.

Our method, also topological and geometric in sprit, is based on
the observation that the computations done by Arai and Mischaikow
in \cite{AM} in the second part of their approach should be in
principle sufficient to obtain the whole result, both the
existence of the homoclinic tangency and its genericity. Obviously
for this end, one should use different tools, ours are essentially
those of differential topology, but developed earlier in the
context of topological dynamics, like the covering relations
combined with the cone conditions, see \cite{KWZ,Z}.  As a test
case we give a computer assisted proof of the existence of the
generic homoclinic tangency for the dissipative H\'enon map with
$b=-0.3$  and $a \approx 1.3145$ in the computation time $0.2$ sec
on the Intel Xeon 5160, 3GHz processor. This should be contrasted
with the fact that in \cite{AM} this case took around 100 minutes.

An application of our method to an ODE requires an efficient
rigorous $C^2$-solver for ODEs. By this we mean an algorithm for
rigorous integration of ODEs together with their  variational
equations  up to the second order. Such an algorithm has been
recently developed by the authors in \cite{WZ2} and is now a part
of the CAPD library \cite{CAPD}. Using this algorithm we were able
to prove the existence of the generic homoclinic tangency for the
$2\pi$-shift along the trajectory of the periodically forced
pendulum equation
\begin{equation}\label{eq:pendulum-intro}
    \ddot{x} + \beta\dot x + \sin(x) = \cos(t).
\end{equation}
for $\beta \approx 0.2471$ (see Theorem~\ref{thm:mainPendulum} in
Section~\ref{sec:pendulum}). The computation time for this proof
is $30$ sec on the Intel Xeon 5160, 3GHz processor.

It should be mentioned that similar results for the H\'enon map have
been obtained by a complex analytic method of Fornaess and Gavosto
\cite{FG1,FG2}. Compared to their method, which depends on the
analyticity of maps, our method and that of Arai and Mischaikow are
rather geometric and topological, and are designed so that they can
be applied to a wider class of maps. Essentially, we require a
continuous family of $C^2$ diffeomorphisms for which we can compute
the image of the maps using interval arithmetic.

The content of the paper may be briefly described as follows: in
Sections~\ref{sec:homoc-proj} and~\ref{sec:quad=transversal}
following mainly~\cite{AM} we give basic definitions and restate
the problem of the existence of the quadratic homoclinic tangency
unfolding generically as the transversality question for a
dynamical system induced by the given one on the projective
bundle. In Section~\ref{sec:covrel} we discuss how this
transversality problem can be solved using the covering relations
linked with the cone conditions. Moreover, we derive computable
estimates for the dependence on a parameter of (un)stable
manifolds of the hyperbolic fixed point, which will be later used
in the computer assisted proofs. In Section~\ref{sec:example}
 we illustrate our approach on a toy
example, with the intention that the reader may see and appreciate
some details of the method, which are later hardly visible when we
report on the computer assisted proof for the H\'enon map and
\textbf{for} the forced damped pendulum in
Sections~\ref{sec:henon} and~\ref{sec:pendulum}, respectively.

\section{Homoclinic tangency and the projectivization}
\label{sec:homoc-proj}

\subsection{Invariant manifold - basic notation}

\begin{definition}
Consider the map $f:X \supset \dom(f) \to X$.

 Let $x \in X$. Any sequence $\{x_k\}_{k \in I}$, where
$I \subset \mathbb{Z}$ is a set containing $0$ and for any  $l_1 <
l_2 < l_3$ in $\mathbb{Z}$ if $l_1,l_3 \in I$, then $l_2 \in I$,
such that
\begin{equation*}
  x_0=x, \qquad f(x_i) =x_{i+1}, \qquad \mbox{for $i,i+1 \in I$}
\end{equation*}
will be called  an orbit through $x$. If $I=\mathbb{Z}_-$, then we
will say that $\{x_k\}_{k \in I}$ is a full backward orbit through
$x$.
\end{definition}

\begin{definition}
Let $X$ be a topological space and  let the map $f:X\supset
\dom(f) \to X$ be continuous.

 Let $Z \subset \mathbb{R}^n$, $x_0 \in Z$, $Z \subset
\dom(f)$. We define
\begin{eqnarray*}
  W^s_Z(z_0,f)&=&\{ z  \: | \: \forall_{n \geq 0} f^n(z)\in Z, \quad \lim_{n\to\infty} f^n(z)=z_0   \}, \\
  W^u_Z(z_0,f)&=&\{ z  \: | \:  \exists \mbox{ $\{x_n\} \subset Z$  a full backward orbit through $z$, such that } \\
     & & \quad \lim_{n\to-\infty} x_n=z_0   \}, \\
  W^s(z_0,f)&=&\{ z  \: | \:  \lim_{n \to \infty} f^n(z)=z_0   \}, \\
 W^u(z_0,f)&=& \{ z  \: | \:  \exists \mbox{ $\{x_n\}$  a full backward orbit through $z$, such that } \\
     & & \quad \lim_{n\to-\infty} x_n=z_0   \}.
\end{eqnarray*}
If $f$ is known from the context, then we will usually drop it and
use $W^s(z_0)$, $W_Z^s(z_0)$ etc.,  instead.
\end{definition}

\subsection{Projectivization of the dynamics}

The notation and setting is the one used in \cite{AM}. Let $f:X
\to X$ be a diffeomorphism of a manifold $X$.
\begin{definition}
\label{def:fphypmap} Let $z_0 \in X$. We say that $z_0$ is a
\emph{hyperbolic fixed point for $f$} iff $f(z_0)=z_0$ and
$\mbox{Sp}(Df(z_0)) \cap S^1=\emptyset$, where $Df(z_0)$ is the
derivative of $f$ at $z_0$, $\mbox{Sp}(A)$ denotes the spectrum of
a  square matrix $A$.
\end{definition}

 We denote the tangent bundle of $X$ by $TX$ and the
differential of $f$ by $Df$. From the dynamical system $f:X \to X$
we can derive a new dynamical system $Pf:PX \to PX$ which is defined
as follows. The space $PX$ is the projective bundle associated to
the tangent bundle of $X$, that is, the fiber bundle on $X$ whose
fiber over $x \in X$ is the projective space of $T_xX$. That is,
\begin{equation*}
  PX=\bigcup_{x \in X}P_xX:=\bigcup_{x\in X}\{\mbox{one-dimensional subspace of
  $T_xX$}\}.
\end{equation*}

For a submanifold  $S$ in $X$ by $PS$ we will denote its
projectivization, which is given by
\begin{equation*}
  PS=\{ (x,[v])\in PX \: | \: x \in S, v \in T_xS \setminus \{0\}  \}
\end{equation*}
It is easy to see that $P(S)$ is a  manifold and
$\mbox{dim}(P(S))=2\mbox{dim}(S)-1$.

Define $Pf$ to be the map induced from $Df$ on $PX$, namely,
$Pf(x,[v]):=(f(x),[Df(x)\cdot v])$  where $0 \neq v \in T_xX$, $[v]$
is the subspace in $T_xX$ spanned by $v$ and $Df(x):T_xX\to
T_{f(x)}X$ is the derivative of $f$ at $x$.

Let us identify $X$ with the zero section of $TX$. We define  the map
 $\pi:TX \setminus X \to PX$ by $\pi((x,v))=(x,[v])$ for $x \in X$, $v \in T_xX$.

Let $p \in X$ be a hyperbolic fixed point of $f$ and let $T_p
X=\tilde{E}^s_p \oplus \tilde{E}^u_p$ be the corresponding
splitting of the tangent space into stable and unstable subspaces
for $Df(p)$.  Define $E^s_p:=\pi({\tilde E}^s_p \setminus \{0\})$
and $E^u_p:=\pi({\tilde E}^u_p \setminus \{0\})$.

\subsection{Dimension two}
\begin{theorem}
\label{thm:Pf-lin} Let $f:\mathbb{R}^2 \to \mathbb{R}^2$ be a
diffeomorphism and let $p$ be a hyperbolic fixed point with one
dimensional stable and unstable manifolds.

Then
\begin{description}
\item[1.]  $(p,E^u_p)$ is a hyperbolic fixed point for $Pf$ such that
\begin{itemize}
\item $W^u((p,E^u_p),Pf)=P(W^u(p,f))$, $\dim W^u((p,E^u_p),Pf)=1$,
\item  $W^s((p,E^u_p),Pf)=\{ (z,[v]) \in PX \: | \:  z \in W^s(p,f), (z,[v]) \notin  P(W^s(p,f))  \}$,
$\dim W^s((p,E^u_p),Pf)=2$.
\end{itemize}
\item[2.]  $(p,E^s_p)$ is a hyperbolic fixed point for $Pf$ such that
\begin{itemize}
\item $W^s((p,E^s_p),Pf)=P(W^s(p,f))$, $\dim W^s((p,E^s_p),Pf)=1$,
\item $W^u((p,E^s_p),Pf)=\{(z,[v]) \in PX \: | \: z \in W^u(p,f), (z,[v]) \notin P(W^u(p,f)) \} $,
 $\dim W^u((p,E^s_p),Pf)=2$.
\end{itemize}
\end{description}

\end{theorem}
\textbf{Proof:} It is easy to see that the stable and unstable
sets are as stated. There remains for us to show the
hyperbolicity, only.  Consider first $(p,E^u_p)$. Let us change
the coordinate system in $\mathbb{R}^2$ such that $p=0$, $Df(p)$
is diagonal
\begin{equation*}
  Df(p)=\begin{pmatrix}
    \lambda & 0 \\
    0 & \mu \
  \end{pmatrix}
\end{equation*}
where $|\lambda| > 1$ and $|\mu| < 1$. From now on all
considerations will be done in this coordinate frame.

Around $(p,E^u(p))$ we will use the coordinate system $\varphi:
\mathbb{R}^2 \times \mathbb{R} \to P(\mathbb{R}^2)$ given by
\begin{equation*}
\varphi(z,v)=(z,[(1,v)]).
\end{equation*}
It is easy to see that
\begin{equation*}
\varphi^{-1}(z,[(v_1,v_2)])=(z,v_2/v_1).
\end{equation*}

In these coordinates $(p,E^u_p)$ is given by $(0,0)$. Now we
compute the linearization of $Pf$ around $(0,0)$.

We have
\begin{eqnarray*}
  Pf(z,v)=(f(z),[Df(z)\cdot(1,v)^T])=(Df(0)\cdot z + o(|z|),
  [Df(z)\cdot(1,v)^T]).
\end{eqnarray*}
Observe that
\begin{multline*}
  Df(z)\cdot(1,v)^T=\begin{pmatrix}
    \frac{\partial f_1}{\partial x}(z) & \frac{\partial f_1}{\partial y}(z) \\
    \frac{\partial f_2}{\partial x}(z) & \frac{\partial f_2}{\partial y}(z) \
  \end{pmatrix}
  \cdot
  \begin{pmatrix}
    1 \\
    v \
  \end{pmatrix}
  =\begin{pmatrix}
    \lambda & 0 \\
    0 & \mu \
  \end{pmatrix}
  \cdot
  \begin{pmatrix}
    1 \\
    v \
  \end{pmatrix}
  + \\
\begin{pmatrix}
    \frac{\partial^2 f_1}{\partial x^2}(0)x  +  \frac{\partial^2 f_1}{\partial x \partial y}(0)y&
    \frac{\partial^2 f_1}{\partial y \partial x}(0)x +  \frac{\partial^2 f_1}{\partial y^2}(0)y \\
    \frac{\partial^2 f_2}{\partial x^2}(0)x +   \frac{\partial^2 f_2}{\partial x \partial y}(0)y &
     \frac{\partial^2 f_2}{\partial y \partial x}(0)x +  \frac{\partial^2 f_2}{\partial y^2}(0)y \
  \end{pmatrix}
  \cdot
  \begin{pmatrix}
    1 \\
    v \
  \end{pmatrix}
  + o(z)  \cdot
  \begin{pmatrix}
    1 \\
    v \
  \end{pmatrix}
  = \\
   \begin{pmatrix}
     \lambda  +  O(|(z,v)|)\\
     \mu v  +   \frac{\partial^2 f_2}{\partial x^2}(0)x +   \frac{\partial^2 f_2}{\partial x \partial y}(0)y
       + o(|(z,v)|) \
   \end{pmatrix}.
\end{multline*}
Now we have to divide the second component of the above vector by
the first one.
Therefore
\begin{eqnarray*}
  v_2/v_1=  \left( \mu v  +   \frac{\partial^2 f_2}{\partial x^2}(0)x +   \frac{\partial^2 f_2}{\partial x \partial y}(0)y
       + o(|(z,v)|)\right)\cdot \lambda^{-1}(1+O(|(z,v)|)= \\
  \frac{\mu}{\lambda} v  +  \frac{1}{\lambda} \frac{\partial^2 f_2}{\partial x^2}(0)x +
 \frac{1}{\lambda}\frac{\partial^2 f_2}{\partial x \partial y}(0)y
 +  o(|(z,v)|).
\end{eqnarray*}
We have proved that the linearization of $Pf$ at $p=(0,0)$ has the
following form
\begin{equation*}
  DPf(p,E^u_p)=\begin{pmatrix}
    \lambda & 0 & 0 \\
    0 & \mu & 0 \\
     \frac{1}{\lambda} \frac{\partial^2 f_2}{\partial x^2}(0) &
     \frac{1}{\lambda}\frac{\partial^2 f_2}{\partial x \partial y}(0) & \frac{\mu}{\lambda} \
  \end{pmatrix}.
\end{equation*}
Hence we see that $(p,E^u_p)$ is a hyperbolic  fixed point for
$Pf$ with one-dimensional unstable and two-dimensional stable
manifolds.

Let us consider now $(p,E^s_p)$. This time we will use the
coordinate system $\varphi: \mathbb{R}^2 \times \mathbb{R} \to
P(\mathbb{R}^2)$ given by
\begin{equation*}
\varphi(z,v)=(z,[(v,1)]).
\end{equation*}
It is easy to see that
\begin{equation*}
\varphi^{-1}(z,[(v_1,v_2)])=(z,v_1/v_2).
\end{equation*}
 Similar computations lead to the following formula for the
linearization of $Pf$ at $(p,E^s_p)$
\begin{equation*}
  DPf(p,E^s_p)=\begin{pmatrix}
    \lambda & 0 & 0 \\
    0 & \mu & 0 \\
     \frac{1}{\mu} \frac{\partial^2 f_1}{\partial x \partial y}(0) &
     \frac{1}{\mu}\frac{\partial^2 f_1}{ \partial y^2}(0) & \frac{\lambda}{\mu} \
  \end{pmatrix}.
\end{equation*}
\qed

From Theorem~\ref{thm:Pf-lin} we obtain the following
\begin{rem}
 If $f$ is as in Theorem~\ref{thm:Pf-lin} and
$W^s(p,f)$ and $W^u(p,f)$ have a nonempty intersection, then
\begin{itemize}
\item if $W^s(p,f)$ and $W^u(p,f)$ are tangent, then we obtain
  a heteroclinic connection from $(p,E^u_p)$ to $(p,E^s_p)$,
\item if $W^s(p,f)$ and $W^u(p,f)$ intersect transversally, then we obtain
  a homoclinic connection from $(p,E^u_p)$ to $(p,E^u_p)$.
\end{itemize}
\end{rem}

\section{Generic unfolding of quadratic tangency as the transversality question}
\label{sec:quad=transversal}

We assume that we have  two  curves in $\mathbb{R}^2$ depending on
some  parameter $a$ and
 given by $u(a,t)=(f_1(a,t),f_2(a,t))$  and
$s(a,t)=(g_1(a,t),g_2(a,t))$. We are interested in establishing
conditions, which will imply the existence of the generic
unfolding of the quadratic tangency between them. Our goal is to
formulate such conditions as the transversality question. This is
Theorem 2.1 from \cite{AM}, where it was stated without proof.

\begin{definition} \cite[Sec.~3.1]{PT}
\label{def:quad-tan} Let $I,J,Z \subset \mathbb{R}$ be intervals.
Let $u_\mu\colon I \to \mathbb{R}^2$ and $s_\mu\colon J \to
\mathbb{R}^2$ for $\mu \in Z $ be two smooth curves depending on
$\mu$ in the smooth way, such that
$u_{\mu_0}(t_u)=s_{\mu_0}(t_s)=q_0$ and $u$ and $s$ are tangent at
$q_0$.

Assume there exists a $\mu$-dependent coordinates in a
neighborhood of $q$ for $\mu$ close to $\mu_0$, such that in these
coordinates we can use $x_1$ (the first coordinate) as the
parameter of our curves and the following holds
\begin{eqnarray}
  s_{\mu}(x_1)&=&(x_1,0), \nonumber\\
  u_{\mu}(x_1)&=&(x_1,ax_1^2 + b(\mu - \mu_0)) \label{eq:form-u-q2}
\end{eqnarray}
where $a \neq 0, b \neq 0$. Then we say that \emph{ the quadratic
tangency of $u$ and $s$ unfolds generically}.
\end{definition}
\begin{rem} It is easy to see that in the above definition we can
exchange the role of curves $u$ and $s$ using the coordinate
transformation given by $\phi_{\mu}(x_1,x_2)=(x_1,x_2 - ax_1^2 -
b(\mu - \mu_0))$.
\end{rem}
\begin{rem}
\label{rem:qtan-2} In the context of Def.~\ref{def:quad-tan}, if
$u_{\mu}(x_1)=g(\mu,x_1)$ then instead of (\ref{eq:form-u-q2}) it
is enough to require
\begin{eqnarray*}
  g(\mu_0,0)=0, \quad \frac{\partial g}{\partial x_1}(\mu_0,0)=0, \\
 \frac{\partial g}{\partial \mu}(\mu_0,0) \neq 0, \quad
  \frac{\partial^2 g}{\partial x_1^2}(\mu_0,0) \neq 0.
\end{eqnarray*}
\end{rem}
\textbf{Proof:} Observe that from our assumptions follows that
\begin{equation*}
  u_\mu(x_1)=(x_1,ax_1^2 + b(\mu-\mu_0) +  \Delta(\mu,x_1)),
\end{equation*}
where
\begin{eqnarray*}
  a\neq 0, \quad b \neq 0, \\
  \Delta(\mu_0,0)=0, \quad   \frac{\partial \Delta}{\partial
  \mu}(\mu_0,0)=0,  \quad \frac{\partial \Delta}{\partial
  x_1}(\mu_0,0)=0.
\end{eqnarray*}
Observe that through the parameter dependent coordinate chart
$\phi_{\mu}(x_1,x_2)=(x_1,x_2 - \Delta(\mu,x_1) )$ we obtain
(\ref{eq:form-u-q2}).
 \qed

\begin{theorem}
\label{thm:trans-quadratic} Let  $\Lambda \subset \mathbb{R}$ be
an interval. Let $u: \Lambda \times \mathbb{R} \to \mathbb{R}^2$
and $s:\Lambda \times \mathbb{R} \to \mathbb{R}^2$ be two
$C^2$-curves depending on parameter $a \in \Lambda$.  Let $a_0 \in
\Lambda$ be a parameter at which curves $u_a$ and $s_a$ are
tangent.

These curves have a quadratic tangency at $a_0$ which unfolds
generically iff there exists $t_u, t_s \in \mathbb{R}$, such that
surfaces
\begin{equation*}
  E(u),E(s):\Lambda \times \mathbb{R} \to \mathbb{R}^3
  \times R\mathbb{P}^1
\end{equation*}
given by
\begin{eqnarray*}
  E(u)(a,t)=\left(a,u(a,t),\left[\frac{\partial u}{\partial t}(a,t)\right]\right)  \\
  E(s)(a,t)=\left(a,s(a,t),\left[\frac{\partial s}{\partial t}(a,t)\right]\right)
\end{eqnarray*}
intersect transversally  at the point
$E(u)(a_0,t_u)=E(s)(a_0,t_s)$.
\end{theorem}
Before the proof we need one lemma.

\begin{lemma}
\label{lem:goodcoord}
  Let $s=(s_1,s_2):\Lambda \times \mathbb{R} \to \mathbb{R}^2$ be a
  $C^2$-map, such that $\frac{\partial s_1}{\partial t}(a_0,t_0) \neq
  0$  for some $(a_0,t_0)$.

  Then there exists a neighborhood  $\Lambda' $ of
  $a_0$, an open set $V$, such that $s(a_0,t_0) \in V$,
  and  $a$-dependent coordinates on $V$ for $a \in \Lambda'$, i.e. $\phi_a:
  V\to  \mathbb{R}^2$ for $a \in \Lambda'$, such that after a suitable
  reparameterization in these new
  coordinates the mapping $s$ has locally the following form
  \begin{equation*}
     s(a,t)=(t,0).
  \end{equation*}
\end{lemma}
\textbf{Proof:} Let us denote by $(x,y)$ coordinates in
$\mathbb{R}^2$. By taking a suitable parameterization and shifting
the coordinates origin and  permuting, if necessary, the
coordinates in $\mathbb{R}^2$ we can assume that $t_0=0$,
$s(a_0,0)=0$ and $\frac{\partial s_1}{\partial t}(a_0,0) > 0$. We
can locally (in a suitable open set $U$) use $x$ as the parameter
of a curve $s(a,\cdot)$ for $a \in \tilde{\Lambda} \subset
\Lambda$. Therefore, we have
\begin{equation*}
  s(a,t)=(t,s_2(a,t)).
\end{equation*}

Consider the following map  $\varphi:\tilde{\Lambda} \times U \to
\tilde{\Lambda} \times \mathbb{R}^2$ given by
\begin{equation*}
  \varphi(a,x,y)=(a,x,y - s_2(a,x)).
\end{equation*}
We have
\begin{equation*}
  D\varphi(a,x,y)=\left[ \begin{array}{ccc}
    1 & 0 & 0 \\
    0 & 1 & 0 \\
    -\frac{\partial s_2}{\partial a}(a,x) & -\frac{\partial s_2}{\partial x}(a,x) & 1 \
  \end{array}
  \right]
\end{equation*}
Therefore, $\varphi$ is a local diffeomorphism. Let open sets
$\Lambda',V$ be such that $\varphi:\Lambda' \times V \to \Lambda'
\times  \mathbb{R}^2 $ is a diffeomorphism on the image and
$(a_0,0) \in \Lambda' \times V$. The $a$-dependent coordinates on
$V$ are given by $\phi_a(x,y)=\varphi(a,x,y)$.

It is easy to see that in the new coordinates given by $\phi_a$
the curve $s$ has the following form
\begin{equation*}
  s(a,t)=(t,0).
\end{equation*}
 \qed

 \noindent
 \textbf{Proof of Theorem~\ref{thm:trans-quadratic}:}
By Lemma~\ref{lem:goodcoord}  we can assume that the curve $s$ is
given by
\begin{equation*}
  s(a,t)=(t,0).
\end{equation*}
Observe that  the transversality implies that in the neighborhood
of the intersection point the curve $u$ can be represented as
follows
\begin{equation*}
  u(a,t)=(t,g_2(a,t)),
\end{equation*}
where $g_2$ satisfies the following conditions
\begin{equation}
g_2(a_0,0)=0, \quad \frac{\partial g_2}{\partial t}(a_0,0)=0.
\label{eq:g2pt0}
\end{equation}

Now we will prove that the transversality of $E(u)$ and $E(s)$  is
equivalent to
\begin{eqnarray*}
  \frac{\partial g_2}{\partial a}(a_0,0) \neq 0, \\
  \frac{\partial^2 g_2}{\partial t^2}(a_0,0) \neq 0.
\end{eqnarray*}
Observe that from Remark~\ref{rem:qtan-2} it follows that the
above conditions together with (\ref{eq:g2pt0}) are equivalent to
the generic unfolding of the quadratic tangency.

Observe that $\frac{\partial s}{\partial t}(a,t)=(1,0)$ and
$\frac{\partial u}{\partial t}(a_0,t=0)=(1,0)$ therefore in the
neighborhood of $[\frac{\partial s}{\partial t}(a,t)]$ we can use
the second coordinate as a chart map in $R\mathbb{P}^1$.

In these coordinates we have
\begin{eqnarray*}
E(s)(a,t)&=&(a,t,0,0)^T, \\
E(u)(a,t)&=&(a,t,g_2(a,t),\frac{\partial g_2}{\partial t}(a,t))^T.
\end{eqnarray*}

We have
\begin{equation*}
  T_{E(s)(a_0,t=0)}=\mbox{span}\begin{pmatrix}
    1 & 0 \\
    0 & 1 \\
    0 & 0 \\
    0 & 0 \
  \end{pmatrix}
\end{equation*}
and
\begin{equation*}
  T_{E(u)(a_0,t=0)}=\mbox{span}\begin{pmatrix}
    1 & 0 \\
    0 & 1 \\
    \frac{\partial g_2}{\partial a}(a_0,0) &  \frac{\partial g_2}{\partial t}(a_0,0) \\
    \frac{\partial^2 g_2}{\partial t\partial a}(a_0,0)  &   \frac{\partial^2 g_2}{\partial t^2}(a_0,0)  \
  \end{pmatrix}.
\end{equation*}
From (\ref{eq:g2pt0}) it follows that the transversality question
is equivalent to  the following determinant being nonzero
\begin{equation*}
\det\begin{pmatrix}
    1 & 0 &    1 & 0 \\
    0 & 1 & 0 & 1\\
    0 & 0 &   \frac{\partial g_2}{\partial a}(a_0,0) &  0 \\
    0 & 0 &  \frac{\partial^2 g_2}{\partial t\partial a}(a_0,0)  &   \frac{\partial^2 g_2}{\partial t^2}(a_0,0)  \
  \end{pmatrix}=  \frac{\partial g_2}{\partial a}(a_0,0)  \cdot \frac{\partial^2 g_2}{\partial
  t^2}(a_0,0).
\end{equation*}
This finishes the proof.
 \qed

\section{How to prove the homoclinic tangency using the covering relations and
 the cone conditions?}

\label{sec:covrel}

We assume that the reader is familiar with the following notions:
h-sets, covering relations, cone conditions and horizontal and
vertical disks as defined in \cite{KWZ,Z}.

We consider a planar map $f_a:\mathbb{R}^2 \supset \dom (f) \to
\mathbb{R}^2$ depending on the parameter $a$, which has a
hyperbolic fixed point $p_a$ with one-dimensional unstable and
stable manifolds. Hence according to the setting from
Section~\ref{sec:homoc-proj} we will work in four-dimensional
space, using coordinates $(a,x,y,v)$, where $a$ is the parameter,
$(x,y) \in \mathbb{R}^2$ and $v$ represents points in
$P(T_{(x,y)}\mathbb{R}^2)$, which we will call the
\emph{tangential} coordinate. We have the map (we abuse the
notation for $Pf$)
\begin{equation*}
  Pf(a,x,y,v)=(a,Pf(x,y,v)).
\end{equation*}

In our method we need the following ingredients
\begin{itemize}
\item the chain of covering relations
\begin{equation}
  N_0 \cover{Pf} N_1 \cover{Pf} \dots \cover{Pf} N_k
  \label{eq:hetero-chain},
\end{equation}
such that the cone conditions are satisfied,
\item $(a,W^u((p_a,E^u_{p_a}), Pf_a))$ as  the horizontal disk in $N_0$ satisfying the cone
conditions,
\item $(a,W^s((p_a,E^s_{p_a}),Pf_a))$ as  the vertical disk in $N_k$ satisfying the cone conditions.
\end{itemize}

If the above conditions are satisfied then from \cite[Theorem
7]{Z} it follows that $(a,W^u((p_a,E^u_{p_a}),Pf_a))$ contains a
horizontal disk satisfying the cone conditions in $N_k$. Hence
$(a,W^u((p_a,E^u_{p_a}),Pf_a))$ and
$(a,W^s((p_a,E^s_{p_a}),Pf_a))$ intersect transversally in $N_k$,
which by Theorem~\ref{thm:trans-quadratic} implies that sets
$W^u((p_a,E^u_{p_a}),Pf_a)$ and $W^s((p_a,E^s_{p_a}),Pf_a)$ have a
quadratic tangency which unfolds generically.

Since the parameter $a$ is not changing under $Pf$, apparently there
is a problem with the realization of the above scenario, because the
covering relations together with the cone conditions imply the
hyperbolicity. The essential point is that this remark is valid for
closed loops of covering relations, but in our setting we just want
a chain of covering relations going from one part of our phase space
to another one. While constructing such chain we will arbitrarily
decide whether we treat the parameter as the "unstable"/"stable"
direction, by simply adjusting the sizes of h-sets in $a$-direction.

To make our scheme to work  we need:
\begin{itemize}
\item to set the dimensions  $u$, $s$ in our h-sets $N_i$ to be equal
to $2$, because this is the dimension of
$(a,W^u((p_a,E^u_{p_a}),Pf_a)$ and $(a,W^s((p_a,E^s_{p_a}),Pf_a)$,
\item $(a,W^u((p_a,E^u_{p_a}),Pf_a)$ should be a horizonal disk in $N_0$, hence
    we should treat $a$ as one of unstable directions,
\item $(a,W^s((p_a,E^s_{p_a}),Pf_a)$ should be a vertical disk in
$N_k$, hence we should treat $a$ as one of stable directions.
\end{itemize}

\subsection{The Lipschitz dependence of stable and unstable manifolds  on parameters  }
\label{subsec:Lipdep}

We need prove that $(a,W^u((p,E^u_p),Pf_a)$ is a horizontal disk
in a suitable h-set $N_0$, this requires at least the Lipschitz
dependence on $a$ of the invariant manifold.  In \cite[Sec.
8.2]{Z} the Lipschitz dependence of (un)stable manifolds with
respect to parameters with explicit and computable constants was
discussed with the eye toward the computer assisted proofs. Here
we will just recall (and refine a bit some estimates) these
results in a form of the cone conditions.

We will be using the norms for quadratic forms (identified in the
sequel with symmetric matrices) which are defined by
\begin{equation*}
  |B(u,v)| \leq \| B \| \| u \| \| v \|.
\end{equation*}
For Euclidian norm  we have
\begin{equation*}
  \| B \| = \max \{|s| \ | \ \mbox{$s$ in an eigenvalue of $B$}\}.
\end{equation*}

\begin{theorem}
\label{thm:coneInvMan} Assume that $(N,Q)$ is an h-set in $\mathbb
R^{u+s}$ with cones and $f_\lambda\colon \mathbb R^{u+s}\to \mathbb
R^{u+s}$ with $\lambda \in C$, where $C$ is a compact interval in
the parameter space and $Q$ has the form  $Q(x,y)= \alpha(x) -
\beta(y)=\sum_{i=1}^u a_i x_i^2 - \sum_{i=1}^sa_{i+u} y_i^2$.

\begin{description}
\item[1.] Assume that for the covering relation $N \cover{f_\lambda} N$  the cone
condition is  satisfied for all $\lambda \in C$.
\item[2.] Let $\epsilon>0$ and $A>0$  be such that for all $\lambda \in
C $ and $z_1,z_2 \in N$ holds
\begin{equation}
  Q(f_\lambda(z_1)-f_\lambda(z_2)) - (1+ \epsilon)Q(z_1 -z_2) \geq
  A(z_1 - z_2)^2.  \label{eq:QlambdaA}
\end{equation}
\item[3.] Let
\begin{eqnarray}
 M&=&\max_{\lambda\in C, z\in N}\left(\sum_i  |a_i|\left\| \frac{\partial \pi_{z_i} f_\lambda}{\partial z}(z) \right\|
  \cdot \left\| \frac{\partial \pi_{z_i}f_\lambda}{\partial \lambda}(z)  \right\|
 \right),  \label{eq:defM} \\
 L&=&\|\beta\|  \cdot \max_{\lambda\in C, z\in N}\left\| \frac{\partial \pi_yf_{\lambda}}{\partial
 \lambda}(z)
    \right\|^2. \label{eq:defL}
\end{eqnarray}
\item[4.] Let $\Gamma >0$  be such that
\begin{equation}
  A - 2M\Gamma - L \Gamma^2 >0.\label{eq:GammaQn}
\end{equation}
\item[5.]
We define
\begin{equation}\label{eq:defDeltaCoeff}
  \delta=\frac{\Gamma^2}{\|\alpha\|}.
\end{equation}
\end{description}

Then  the set $W^s_N(p_\lambda,f_\lambda)$ for $\lambda \in C$ can
be parameterized as a vertical disk in $C \times N$ for the
quadratic form $\tilde{Q}(\lambda,z)=\delta Q(z) - \lambda^2$.
\end{theorem}

Before the proof let us make two observations concerning constants
$A,\epsilon,\Gamma$.

\begin{rem}
The existence of  $A$ and $\epsilon$ in  (\ref{eq:QlambdaA}) is a
consequence of the cone condition for covering relation $N
\cover{f_\lambda} N$.
 We would like to have as big $A$ as possible. This forces $\epsilon \to 0$,
but $\epsilon$ is not used in the sequel.
\end{rem}

\begin{rem}
Since $A>0$, therefore $\Gamma$ in (\ref{eq:GammaQn}) always
exists, but it is desirable to look for largest $\Gamma$ possible.
\end{rem}

\noindent \textbf{Proof of Theorem~\ref{thm:coneInvMan}:} We would
like to obtain that for $|\lambda_1 - \lambda_2| \leq \Gamma \|z_1
- z_2 \|$ holds
\begin{equation*}
   Q(f_{\lambda_1}(z_1) - f_{\lambda_2}(z_2) ) >
   (1+\epsilon)Q(z_1-z_2).
\end{equation*}

Let $B$ be a unique symmetric form, such that $B(u,u)=Q(u)$. Observe
that
\begin{eqnarray*}
  Q(f_{\lambda_1}(z_1) - f_{\lambda_2}(z_2) ) - (1+\epsilon)Q(z_1-z_2)
  =\\
   Q(f_{\lambda_1}(z_1) - f_{\lambda_1}(z_2) ) -
 (1+\epsilon) Q(z_1-z_2) + \\
   2B(f_{\lambda_1}(z_1) - f_{\lambda_1}(z_2), f_{\lambda_1}(z_2) -
   f_{\lambda_2}(z_2)) + Q(f_{\lambda_1}(z_2) -
   f_{\lambda_2}(z_2)).
\end{eqnarray*}
The first term in the above expression will be estimated using
(\ref{eq:QlambdaA}).

For the third term we obtain
\begin{eqnarray*}
  Q(f_{\lambda_1}(z_2)-f_{\lambda_2}(z_2)) \geq -
  \beta(\pi_y(f_{\lambda_1}(z_2)-f_{\lambda_2}(z_2))) \geq \\
    -\|\beta\| \cdot \max_{\lambda \in C}\left\| \frac{\partial \pi_yf_{\lambda}}{\partial
    \lambda}(z_2)
    \right\|^2 \cdot (\lambda_1 - \lambda_2)^2 \geq\\
     -\|\beta\| \max_{\lambda\in C, z\in N} \left\| \frac{\partial \pi_yf_{\lambda}}{\partial
     \lambda}(z)
    \right\|^2 \cdot \Gamma^2 \|z_1 - z_2\|^2 = -L\Gamma^2 \|z_1-z_2\|^2.
\end{eqnarray*}
Finally, for  the second term we have
\begin{eqnarray*}
 | B(f_{\lambda_1}(z_1) -
  f_{\lambda_1}(z_2),f_{\lambda_1}(z_2)-f_{\lambda_2}(z_2)) | \leq
  \\
  \max_{\lambda\in C, z\in N} \left(\sum_i  |a_i|\left\| \frac{\partial \pi_{z_i} f_\lambda}{\partial z}(z) \right\|
  \cdot \left\| \frac{\partial \pi_{z_i}f_\lambda}{\partial \lambda}(z)  \right\|
 \right) \cdot \Gamma \| z_1 - z_2 \|^2.
\end{eqnarray*}

From the above computations and (\ref{eq:defM}--\ref{eq:defL}) we
obtain the following
\begin{equation*}
  Q(f_{\lambda_1}(z_1) - f_{\lambda_2}(z_2) ) -
  (1+\epsilon)Q(z_1-z_2) \geq  \left(A - 2M\Gamma - L
  \Gamma^2\right)\|z_1-z_2\|^2.
\end{equation*}

For $\Gamma$ and $\delta$ as in  (\ref{eq:GammaQn}) and
(\ref{eq:defDeltaCoeff}) it is proved in \cite[ Lemma 23]{Z}
that, for  $\lambda_1,\lambda_2\in C$, $\lambda_1 \neq \lambda_2$
and $z_i \in W^s_N(p_{\lambda_i},f_{\lambda_i})$ for $i=1,2$ holds
\begin{equation*}
  \delta Q(z_1 - z_2) - (\lambda_1 - \lambda_2)^2 < 0.
\end{equation*}

Hence $W_N^s(p_\lambda,f_\lambda)$ for $\lambda \in C$ is vertical
disk in $C \times N$ for the quadratic form
$\tilde{Q}(\lambda,z)=\delta Q(z) - \lambda^2$.
 \qed

\textbf{Comments:}
\begin{itemize}
\item Since $\tilde{Q}(x,y,\lambda)= \delta \alpha(x) - \delta \beta(y) -
  \lambda^2$, hence due to the negative sign in front of $\lambda^2$ it follows that
   when trying to represent  the stable manifold as the vertical disk in $C \times N$,
    we must treat
   the parameter  as  a 'stable' direction in an h-set.
\item For the unstable manifold  we have to take the inverse map and
we obtain different $\delta$. Since taking the inverse involves
changing the sign of $Q$ we end up with the following quadratic
form
\begin{equation*}
\tilde{Q}(x,y,\lambda)=\delta \alpha(x) + \lambda^2 - \delta
\beta(y).
\end{equation*}
 Looking at sign in front of $\lambda^2$ we see that $\lambda$ appears as an 'unstable' direction
 in $C \times N$.
\end{itemize}

\section{A toy example }
\label{sec:example}  In this section we show how to construct the
chain of covering relations (\ref{eq:hetero-chain}) discussed in
the first part of Section~\ref{sec:covrel} for a special model map
with a quadratic tangency which unfolds generically. Our intention
is  that the reader may see and appreciate some details of the
method, which are later hardly visible when we report on the
computer assisted proofs for the H\'enon map and the forced
pendulum in Sections~\ref{sec:henon} and~\ref{sec:pendulum},
respectively.

 We define a map $f:\mathbb{R}^2
\to \mathbb{R}^2$ depending on the parameter $a$ as follows:
\begin{itemize}
\item $(0,0)$ is a hyperbolic fixed point and in a  neighborhood
   of $(0,0)$ the map $f_a$ is linear
   \begin{equation*}
      f_a(x,y)=(\lambda x,\mu y)
   \end{equation*}
   where $|\lambda| >1$ and $|\mu| < 1$,
\item in a  neighborhood of the point $(1,0)$ we have the homoclinic tangency
  for $a=0$. We assume that $f_a$ acts as follows
     \begin{equation*}
       f_a(1+x,y)=(x^2+y+a,1-x).
     \end{equation*}
\end{itemize}
When compared with the full problem of establishing the existence
of generic unfolding of homoclinic tangency for the H\'enon map
and the forced damped pendulum, the above model map avoids the
problems related to providing the explicit estimates on the
dependence of $W^{s,u}$ on the parameter, which are given in
Theorem~\ref{thm:coneInvMan}.

\subsection{Chain of covering relations }
Let $p=(0,0)$, then $E^u_p=[(1,0)]$ and $E^s_p=[(0,1)]$. Our goal
is to construct a chain of covering relations 'linking'
$(p,E^u_p)$ with $(p,E^s_p)$.

\textbf{The beginning of the chain:}

 Let us first see how the
covering relations look in a  neighborhood of $(p,E^u_p)$. We use
the chart $(x,y,v,a)$, where $v \mapsto [(1,v)]$. Therefore the
map $F=(Pf_a,a)$  works as follows
\begin{equation*}
  F(x,y,v,a)=(\lambda x,\mu y,(\mu/\lambda) v,a).
\end{equation*}

We define an $h$-set, so that $(x,a)$ are the 'unstable' directions,
$(y,v)$ are the 'stable' ones. It is easy to see, that there exists
a sequence of $h$-sets $N_i=(c_i +[-x_i,x_i]) \times [-y_i,y_i]
\times [-v_i,v_i] \times [-a_i,a_i]$, with  $c_0=0$, $c_k=1$ and
such that
\begin{equation*}
   N_0 \cover{F} N_1 \cover{F} \dots \cover{F} N_k.
\end{equation*}
Observe that the necessary conditions are
\begin{eqnarray*}
     c_{i+1}+[-x_{i+1},x_{i+1}] &\subset& \lambda c_i + [-|\lambda| x_i,
     |\lambda| x_i], \\
     |\mu| y_i &<& y_{i+1}, \\
     \left| \mu/\lambda  \right|  v_i &<& v_{i+1}, \\
     a_i &>& a_{i+1}.
\end{eqnarray*}

\textbf{The end of the chain}:

 In a  neighborhood of $(p,E^s_p)$ the 'unstable' subspace is $(x,w)$,
 where $w$ corresponds to the $[(w,1)]$. The map $F$ works as follows
\begin{equation*}
   F(x,y,w,a)=(\lambda x,\mu y,(\lambda/\mu) w,a).
\end{equation*}
It is easy to see, that there exists a sequence of $h$-sets
$M_i=(\bar{c}_i +[-\bar{x}_i,\bar{x}_i]) \times
[-\bar{y}_i,\bar{y}_i] \times [-\bar{w}_i,\bar{w}_i] \times
[-\bar{a}_i,\bar{a}_i]$, such that
\begin{equation*}
  \bar{c}_0=0, \  \bar{c}_s=1, \qquad  M_s \cover{F} M_{s-1} \cover{F} \dots  \cover{F}
  M_0.
\end{equation*}
Observe that the necessary conditions are
\begin{eqnarray*}
 (\bar{c}_{i+1} +[-\bar{x}_{i+1},\bar{x}_{i+1}]) &\subset&  \mu\bar{c}_i
 +[-|\mu|\bar{x}_i,|\mu|\bar{x}_i], \\
 \bar{y}_{i+1} &>& |\mu| \bar{y}_i, \\
 \bar{w}_{i+1} &<& \left| \lambda/\mu \right| \bar{w}_i, \\
 \bar{a}_{i+1} &<& \bar{a}_i.
\end{eqnarray*}

\textbf{Switching from unstable to the stable manifold.} We want
the covering relation $N_k \cover{F} M_s$, where $N_k$ and $M_s$
are as above.

In $N_k$ the nominally unstable directions are $(x,a)$, while  in
$M_s$  parameterized by $(x,y,w,a)$, where $w \mapsto [(w,1)]$ the
'unstable' direction is $(x,w)$. We have
\begin{equation*}
  F(1+x,y,v,a)=(x^2+y+a,1-x,-(2x+v),a).
\end{equation*}
We would like to have the following  homotopy for covering
relation
\begin{equation*}
  G_t(1+x,y,v,a)=(a+t(x^2 + y), 1 - tx, -2x -tv, ta).
\end{equation*}

Therefore if we want a covering relation $N_k \cover{F} M_s$,
where $N_k=(1,0,0,0)+ [-x_k,x_k]\times
[-y_k,y_k]\times[-v_k,v_k]\times [-a_k,a_k]$ and $M_s=(0,1,0,0)+
[-\bar{x}_s,\bar{x}_s]\times
[-\bar{y}_s,\bar{y}_s]\times[-\bar{w}_s,\bar{w}_k]\times
[-\bar{a}_k,\bar{a}_k]$, then we need to satisfy the following set
of inequalities
\begin{eqnarray*}
  a_k - x_k^2 - y_k > \bar{x}_s, \\
  2x_k - v_k > \bar{w}_s, \\
  \bar{y}_s > x_k, \\
  \bar{a}_s > a_k.
\end{eqnarray*}
It is easy to see that this set of inequalities has a solution.
For example if $a_k=\Delta < 1$, then we can choose (for some
small $\epsilon >0$)
\begin{eqnarray*}
  x_k= \Delta/2, \\
  y_k=\bar{x}_s = \Delta/3, \\
  v_k=\bar{w}_s=(1-\epsilon)\Delta/2, \\
  \bar{y}_s= (0.5+\epsilon) \Delta, \\
  \bar{a}_s=(1+\epsilon) \Delta.
\end{eqnarray*}
Observe that the expanding direction $x$ is stretched in
$w$-direction and $a$ direction is stretched across $x$-direction
(in the target set).

It is clear that we can easily build the desired chain of covering
relation of the form
\begin{equation*}
  N_0 \cover{F} N_1 \cover{F} \dots \cover{F} N_k \cover{F}   M_s \cover{F} M_{s-1} \cover{F} \dots  \cover{F}
  M_0.
\end{equation*}

\subsection{Cone conditions}

It turns out  that the cone conditions at the beginning and the
end of the chain (\ref{eq:hetero-chain})  are relatively easy to
satisfy in the situation when dynamics is linear (there is no
parameter dependence) in a  neighborhood of the hyperbolic fixed
point, otherwise the issue becomes delicate see
Theorem~\ref{thm:coneInvMan} in Section~\ref{subsec:Lipdep}.

As a rule in this subsection the coordinate order for h-sets $N,M$
will such that always the nominally unstable coordinates are
written first.

\textbf{At the beginning of the chain} For covering relations $N_i
\cover{F} N_{i+1}$ the expanding directions are $(x,a)$. Let $N_i$
be an h-set with $Q_i$-cones given by
\begin{equation*}
  Q_{N_i}(x,a,y,v)=\alpha_i x^2 + \beta_i a^2 - \gamma_i y^2  -
  \delta_i v^2.
\end{equation*}
The map $F$ is linear for $N_i \cover{F} N_{i+1}$ therefore it is
enough to check whether
\begin{equation*}
  Q_{N_{i+1}}(F(x,a,y,v)) > Q_{N_i}(x,a,y,v).
\end{equation*}
We have
\begin{eqnarray*}
  Q_{N_{i+1}}(F(x,a,y,v)) - Q_{N_i}(x,a,y,v) = \\ \alpha_{i+1} \lambda^2  x^2 + \beta_{i+1} a^2 -
   \gamma_{i+1} \mu^2 y^2  -
   (\mu/\lambda)^2 \delta_{i+1} v^2  - \\
   (\alpha_i x^2 + \beta_i a^2 - \gamma_i y^2  -
  \delta_i v^2)= \\
   (\alpha_{i+1}\lambda^2 - \alpha_i)x^2 + (\beta_{i+1} -
   \beta_i)a^2 + (\gamma_i - \mu^2 \gamma_{i+1})y^2 + (\delta_i - (\mu/\lambda)^2 \delta_{i+1})
\end{eqnarray*}
Hence we will have cone conditions satisfied, when for example
$\alpha_{i+1}=\alpha_i$, $\gamma_{i+1}=\gamma_i$,
$\delta_{i+1}=\delta_i$. The only strict requirement is
\begin{equation*}
  \beta_{i+1} > \beta_i.
\end{equation*}

\textbf{At the end of the chain.} It goes through as easily as the
previous case.  For covering relations $M_{i} \cover{F} M_{i-1}$ the
unstable directions are $(x,w)$. We set
\begin{equation*}
  Q_{M_i}(x,w,y,a)=A_ix^2 +B_i w^2  - C_iy^2 -D_i a^2.
\end{equation*}
As above it is enough to whether
\begin{equation*}
  Q_{M_{i-1}}(F(x,w,y,a)) > Q_{M_i}(x,w,y,a).
\end{equation*}
We have
\begin{multline*}
  Q_{M_{i-1}}(F(x,w,y,a)) - Q_{M_i}(x,w,y,a)= \\
  (A_{i-1}\lambda^2 -
  A_i)x^2 + (B_{i-1}(\lambda/\mu)^2 - 1 B_i)w^2 +\\
   (C_{i-1} - C_i
  \mu^2)y^2 + (D_{i-1} - D_i)a^2.
\end{multline*}
We see that we need to have
\begin{equation*}
  D_{i-1} > D_i.
\end{equation*}
For the remaining coefficients we can set
\begin{equation*}
   A_i=A_{i-1}, B_{i-1}=B_i, C_{i-1}=C_i.
\end{equation*}

\textbf{Switching from the unstable to the stable manifold.}
Consider the covering relation  $N_k \cover{F} M_s$. This time we
are in the nonlinear regime. We use on $N_k$ coordinates
$(x,a,y,v)$ and on $M_s$ coordinates $(x,w,y,a)$

This means that in these coordinates (we remove the shifts of the
origin of coordinate frame to $(1,0,0,0)$ and $(0,0,1,0)$ for
$N_k$ and $M_s$ respectively)
\begin{equation*}
  F(x,a,y,v)=(x^2+y+a,-2x - v,-x,a).
\end{equation*}

We set
\begin{eqnarray*}
  Q_{N}(x,a,y,v)&=&\alpha x^2 + \beta a^2 - \gamma y^2  -
  \delta v^2, \\
  Q_{M}(x,w,y,a)&=&Ax^2 +B w^2  - Cy^2 -D a^2.
\end{eqnarray*}

We have
\begin{equation*}
  DF(x,a,y,v)=\begin{bmatrix}
    2x & 1 & 1 & 0 \\
    -2 & 0 & 0 & -1 \\
    -1 & 0 & 0 & 0 \\
    0 & 1 & 0 & 0  \
  \end{bmatrix}.
\end{equation*}
After some computations we  obtain
\begin{eqnarray*}
 \tilde{Q}= (DF)^T Q_M DF - Q_N =
  \begin{bmatrix}
    4Ax^2 + (4B - C - \alpha)& 2Ax & 2Ax & 2B \\
    2Ax & A-D - \beta & A & 0 \\
    2Ax & A & A + \gamma & 0 \\
    2B & 0 & 0 & B + \delta  \
  \end{bmatrix}.
\end{eqnarray*}
We are interested whether $\tilde{Q}$ is positive definite. Observe
that since the positive definiteness is an open condition we will
show that $\tilde{Q}$ is positive definite for $x=0$ and then we
will know that the same holds for $|x|$ small.

After we set $x=0$ we obtain
\begin{eqnarray*}
 \tilde{Q}=
  \begin{bmatrix}
    4B - C - \alpha& 0 & 0 & 2B \\
    0 & A-D - \beta & A & 0 \\
    0 & A & A + \gamma & 0 \\
    2B & 0 & 0 & B + \delta  \
  \end{bmatrix}.
\end{eqnarray*}
It is easy to see that after rearrangement of coordinates the
question is reduced to the positive definiteness  of the following
two matrices
\begin{eqnarray*}
  \tilde{Q}_1=   \begin{bmatrix}
    4B - C - \alpha&  2B \\
    2B &  B + \delta  \
  \end{bmatrix},
  \quad \tilde{Q}_2=
  \begin{bmatrix}
    A-D - \beta & A  \\
     A & A + \gamma  \
  \end{bmatrix}.
\end{eqnarray*}
For example, we can  set $A=B=C=1$ and $D=\frac{1}{2}$, then we
obtain
\begin{eqnarray*}
  \tilde{Q}_1=   \begin{bmatrix}
    3 - \alpha&  2 \\
    2 &  1 + \delta  \
  \end{bmatrix},
  \quad \tilde{Q}_2=
  \begin{bmatrix}
    1/2 - \beta & 1  \\
     1 & 1 + \gamma  \
  \end{bmatrix}.
\end{eqnarray*}
It is now easy to see that we get what we want, when we set
\begin{eqnarray*}
  \alpha=1, \quad \delta > 1, \\
  \beta=1/4, \quad \gamma > 3.
\end{eqnarray*}
Summarizing we obtained
\begin{eqnarray*}
  Q_{N}(x,a,y,v)&=& x^2 + a^2/4 - 4 y^2  -
  2 v^2, \\
  Q_{M}(x,w,y,a)&=&x^2 + w^2  - y^2 -  a^2/2.
\end{eqnarray*}

\section{Application to the H\'enon map.}

\label{sec:henon}

In this section we will show that the method introduced in the
previous sections can be successfully applied to a specific system.
Let us consider the H\'enon map
\begin{equation}\label{eq:henon}
H_{a,b}(x,y) = (a-x^2 +by,x).
\end{equation}
The following theorems have been proven in \cite{AM}.
\begin{theorem}{\cite[Thm.1.1]{AM}}
\label{thm:henon-q-b} There exists an open neighborhood $B$ of
parameter value $b=0.3$ such that for each parameter $b\in B$
there exists a parameter value
\begin{equation*}
a\in[1.392419807915, 1.392419807931]
\end{equation*}
such that the H\'enon map has a quadratic tangency unfolding
generically for the fixed point
\begin{equation*}
x_{a,b}=y_{a,b} = -\frac{1}{2}\left(b-\sqrt{(b-1)^2+4a}-1\right).
\end{equation*}
\end{theorem}
\begin{theorem}{\cite[Thm.1.2]{AM}}
There exists an open neighborhood $B$ of parameter value $b=-0.3$
such that for each parameter $b\in B$ there exists a parameter value
\begin{equation}\label{eq:AMParameter}
a\in[1.314527109319, 1.314527109334]
\end{equation}
such that the H\'enon map has a quadratic tangency unfolding
generically for the fixed point
\begin{equation*}
x_{a,b}=y_{a,b} = \frac{1}{2}\left(b-\sqrt{(b-1)^2+4a}-1\right).
\end{equation*}
\end{theorem}

\begin{figure}[htbp]
\begin{center}
\includegraphics[width=.5\textwidth]{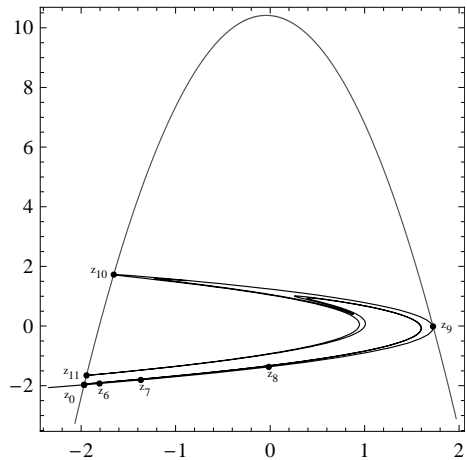}
\end{center}
\caption{Parts of the unstable and stable manifolds of the H\'enon
map at the fixed point for $b=-0.3$ and
$a=a_0=1.3145271093265$.\label{fig:henonManifolds}}
\end{figure}

The numerical evidence of the existence of a homoclinic tangency for
some parameter values $(a,b)\approx(1.3145271093265,-0.3)$ is shown
on Fig.~\ref{fig:henonManifolds} -- see also \cite[Fig.~1.1]{AM}

The main motivation for us to study the existence of homoclinic
tangencies for the H\'enon map was to verify if our method could
work in this relatively easy example.

Denote by
\begin{equation*}
a_0 = 1.3145271093265
\end{equation*}
a center of the interval (\ref{eq:AMParameter}). The aim of this
section is to prove the following theorem
\begin{theorem}\label{thm:mainHenon}
There exists an open neighborhood $B$ of the parameter value
$b=-0.3$ such that for each $b\in B$ there is a parameter $a\in
a_0+[-10^{-5},10^{-5}]$ such that the H\'enon map \eqref{eq:henon}
$H_{a,b}$ has a quadratic homoclinic tangency unfolding generically
for the fixed point
\begin{equation}\label{eq:fixPoint}
x_{a,b}=y_{a,b} =
\frac{1}{2}\left(b-\sqrt{(b-1)^2+4a}-1\right)\approx
-1.9679632427827796.
\end{equation}
\end{theorem}

The authors in \cite{AM} report relatively long computational time
(approximately 100 minutes) for the computer assisted proof of
\cite[Thm.1.2]{AM}. It turns out that using our method the
verification of necessary inequalities in the computer assisted
proof of Theorem~\ref{thm:mainHenon} has been completed in $0.2$
second on the Intel Xeon 5160, 3GHz processor. This very good
efficiency allowed us to apply the method to a map coming from an
ODE - a suitable Poincar\'e map in the forced-damped pendulum. The
details will be given in the next section.

Computer assisted proof of Theorem~\ref{thm:mainHenon} will be
presented in the following subsections in which we verify:
\begin{itemize}
\item the existence of a heteroclinic chain of covering relations for $PH$,
\item  the cone conditions  along the above mentioned
heteroclinic chain,
\item the  cone conditions at the beginning and at the
end of the chain which allow us to parameterize the center-unstable
and center-stable manifolds as a horizontal or vertical discs,
respectively, in proper h-sets.
\end{itemize}

\subsection{The existence of a heteroclinic chain of covering relations for $PH$.}
In Theorem~\ref{thm:mainHenon} we chose the center of the interval
(\ref{eq:AMParameter}) $a_0=1.3145271093265$ as a good candidate for
the homoclinic tangency parameter corresponding to $b_0=-0.3$.

In order to define the sets which will appear in the heteroclinic
chain, we need to set a local chart for the manifold $P\mathbb
R^2$. It turns out, that it is enough for our purpose to use the
parametrization
\begin{equation}\label{eq:atlas}
    \psi:\mathbb R^2\times (0,\pi)
    \ni(x,y,t)\mapsto (x,y,[(\cos(t),\sin(t))])\in
    P\mathbb R^2.
\end{equation}
This parametrization excludes one point on the manifold, but
trough our computations this point does not appear either as an
argument or a value of $PH_a$. All the sets will be expressed
using coordinates  $(x,y,t,a)$.

Put
\begin{equation}\label{eq:extendedHenon}
    PH(x,y,t,a) = \left((\psi^{-1}\circ PH_{a,b_0}\circ \psi)(x,y,t),a\right)
\end{equation}
In the sequel by $\pi_t$ we will denote a projection onto tangent
coordinate, i.e.
\begin{equation*}
\pi_t(x,y,[u]) = \frac{u}{\|u\|}.
\end{equation*}
With some abuse of notation we will use the same symbols for the
projections
\begin{equation*}
\pi_t(x,y,[u],a) = \frac{u}{\|u\|}\quad\text{or}\quad
\pi_t(x,y,t,a)=(\cos(t),\sin(t))
\end{equation*}
but it is clear from the list of arguments which projection has to
be used. In each case the value of $\pi_t$ is a vector.

To simplify the notation we will use $z_0=(x_0,y_0) =
(x_{a_0,b_0},y_{a_0,b_0})$ as defined in \eqref{eq:fixPoint}.
Since we always will have the fixed value of the parameter
$b_0=-0.3$ we will write $H_a$ instead of $H_{a,b_0}$ in the
sequel. Let $u_0$ and $s_0$ be normalized with respect to the
Euclidean norm eigenvectors of $DH_{a_0}(z_0)$ given explicitly by
\begin{equation}\label{eq:HenonStableUnstable}
\begin{array}{l}
 u_0 =
\frac{(-x_0+\sqrt{x_0^2+b_0},1)}{\left\|(-x_0+\sqrt{x_0^2+b_0},1)\right\|}
\approx (0.9680131177714217873,0.250899589123719882),\\
s_0 =
\frac{(-x_0-\sqrt{x_0^2+b_0},1)}{\left\|(-x_0-\sqrt{x_0^2+b_0},1)\right\|}
\approx -(0.07752307795993337433,0.996990557820693689)
\end{array}
\end{equation}
and let $M=[u_0^T,s_0^T]$ be a matrix of eigenvectors. Put
\begin{equation*}
    z_1 = z_0 + 0.0001993152279412426 u_0 + 2.50404\cdot
    10^{-11}s_0.
\end{equation*}
The above point has been chosen as a good approximation of the
homoclinic tangency point for $H_{a_0}$. Namely, we have
\begin{eqnarray}
    \|H_{a_0}^{-1}(z_1) - z_0\|\leq 5.2\cdot 10^{-5},\nonumber\\
    \|H_{a_0}^{14}(z_1) - z_0\|\leq 1.2\cdot 10^{-5},\nonumber\\
    M^{-1}\left(\pi_t\left(PH_{a_0}^{14}(z_1,[u_0])\right)\right)\approx
     (-4.71\cdot10^{-7},0.999999847)\label{eq:iniPointAccuracy}.
\end{eqnarray}

We see that the unstable direction $u_0$ at $z_1$ is mapped under
$14$-th iterate  of $H_{a_0}$ very close to the stable direction
$s_0$ at the point $H^{14}_{a_0}(z_1)$. Let us underline that to
get the estimation (\ref{eq:iniPointAccuracy}) we need to compute
$DH_{a_0}^{14}$ with at least {\tt long double} precision of the
floating point arithmetics.

The points on the trajectory of $z_1$ will be the centers of the
sets which appear in the computer assisted proof. Put
\begin{eqnarray*}
    c_0 &=& (\psi^{-1}(z_0,u_0),a_0),\\
    c_1 &=& (\psi^{-1}(z_1,u_0),a_0),\\
    c_{i+1} &=& PH^i(c_1), \quad\text{for } i=1,\ldots,13,\\
    c_{15} &=& (\psi^{-1}(z_0,s_0),a_0).
\end{eqnarray*}
For further use we set also
\begin{eqnarray*}
    z_{i+1}&=& H_{a_0}^i(z_1), \quad\text{for } 1=1,\ldots,13,\\
    z_{15} &=& z_0.
\end{eqnarray*}
Some of these points are shown in Fig.~\ref{fig:henonManifolds}.

 Now, we have to chose an approximate
stable and unstable directions at $c_i$, $i=0,\ldots,15$. On each
set centered at $c_i$, the coordinate system will be given by a
matrix
\begin{equation}\label{eq:HenonMatrices}
    M_i = \begin{bmatrix}
        (u_i)_1 & (s_i)_1 & 0 & 0\\
        (u_i)_2 & (s_i)_2 & 0 & 0\\
        0 & 0 & 1 & 0\\
        0 & 0 & 0 & 1
        \end{bmatrix}
\end{equation}
where $u_0$, $s_0$ are given by (\ref{eq:HenonStableUnstable}) and
$u_i$ and $s_i$ are computed as follows.
\begin{itemize}
\item Put $u_{15}=u_0$, $s_{15}=s_0$, $u_1=u_0$, $s_1=s_0$.
\item For $i=2,\ldots,8$ we set
\begin{eqnarray}
    u_i&=&\pi_t(c_i),\label{eq:unstFirstChain}\\
    s_i&=&\pi_t(PH^{-1}(z_{i+1},\pi_t(c_{i+1})^\perp)).
\end{eqnarray}
From numerical simulations we get, that between the points $c_8$ and
$c_9$ the role of the tangent coordinate is changing from
contracting to expanding. Therefore, the unstable direction
propagates very well on these sets just by (\ref{eq:unstFirstChain})
for $i=2,\ldots,8$. For the inverse map, the unstable directions
$u_i$ become repelling. Therefore the preimage of an orthogonal
direction to $u_{i+1}$ becomes a good enough approximation for our
method of the stable direction at $z_i$.
\item According to a good choice of the homoclinic point (\ref{eq:iniPointAccuracy}),
the tangent coordinate at $c_i$ become a good
approximation of the stable direction on sets centered at
$c_9,\ldots,c_{14}$. Therefore, we can set
\begin{eqnarray*}
    s_{i} &=& \pi_t(c_i),\quad\text{for } i=9,\ldots,14,\\
    u_{9} &=& \pi_t(PH_{a_0}(z_8,[s_8])),\\
    u_{i+1} &=& \pi_t(PH_{a_0}(z_i,[u_i])),\quad\text{for }i=9,\ldots,13.
\end{eqnarray*}
\end{itemize}

There remains for us to set the sizes of the sets. We define the
h-sets by
\begin{equation}
    N_i = c_i + M_i\cdot\left(d_i\cdot[-1,1]\right),\label{eq:defHSetsHenon}
\end{equation}
where the diameters $d_i$ are listed in Table~\ref{tab:HenonHSets}
and chosen from numerical experiments as well working. On h-sets
$N_0,\ldots,N_8$ the expanding directions are $u$ and $a$
coordinates, while on $N_9,\ldots,N_{15}$ the expanding directions
are $u$ and $t$ coordinates.
\begin{table}
\begin{center}
\begin{tabular}{|c|c|c|c|c|}
    \hline
        $i$ & $10^5\cdot (d_i)_1$& $10^5\cdot (d_i)_2$ & $10^5\cdot (d_i)_3$ & $10^5\cdot (d_i)_4$\\
            & unstable dir. & stable dir. & tangent dir. & parameter\\
    \hline
        $0$ & $7$ & $1$ & $2$ & $(1.01)^8$\\
    \hline
        $1$ & $1$ & $1$ & $2$ & $(1.01)^7$\\
    \hline
        $2$ & $1$ & $1$ & $2$ & $(1.01)^6$\\
    \hline
        $3$ & $1$ & $1$ & $2$ & $(1.01)^5$\\
    \hline
        $4$ & $1$ & $1$ & $2$ & $(1.01)^4$\\
    \hline
        $5$ & $1$ & $1$ & $2$ & $(1.01)^3$\\
    \hline
        $6$ & $1$ & $1$ & $2$ & $(1.01)^2$\\
    \hline
        $7$ & $1$ & $1$ & $2$ & $1.01$\\
    \hline
        $8$ & $1$ & $1$ & $2$ & $1$\\
    \hline
        $9$ & $0.5$ & $1.25$ & $0.25$ & $1.01$\\
    \hline
        $10$ & $0.75$ & $1.25$ & $0.25$ & $(1.01)^2$\\
    \hline
        $11$ & $1$ & $1.25$ & $0.25$ & $(1.01)^3$\\
    \hline
        $12$ & $1$ & $1.25$ & $0.25$ & $(1.01)^4$\\
    \hline
        $13$ & $1$ & $1.25$ & $0.25$ & $(1.01)^5$\\
    \hline
        $14$ & $1$ & $1.25$ & $0.25$ & $(1.01)^6$\\
    \hline
        $15$ & $1$ & $2$ & $0.25$ & $(1.01)^7$\\
    \hline
\end{tabular}
\end{center}
\caption{Diameters of the h-sets in the heteroclinic sequence for
the H\'enon map. The diameters in the table are scaled by the factor
$10^{5}$. On h-sets $N_0,\ldots,N_8$ the expanding directions are
$u$ and $a$ coordinates, while on $N_9,\ldots,N_{15}$ the expanding
directions are $u$ and $t$ coordinates. \label{tab:HenonHSets}}
\end{table}

Let us comment briefly about the choice of the sizes of the sets
presented in Table~\ref{tab:HenonHSets}. For the sets
$N_0,\ldots,N_8$ the tangent coordinate is chosen as a contracting
direction. Therefore the set in this direction must be large enough
(here $2\cdot10^{-5}$) to be able to easy verify that we have
contraction on this variable. On the other hand, for the sets
$N_9,\ldots,N_{15}$ we have a strong expansion on this coordinate.
Therefore the diameter of these sets on this coordinate is smaller
here because it is easier to verify that the image of some walls of
the previous set in the chain is outside of a small set than inside
of it.

Moreover, the tangent coordinate of $N_9$ has to be covered by the
unstable coordinate of $N_8$. This expansion is weak and forces
decreasing of the size $N_9$ on the tangent coordinate.

The unstable and stable direction are comparable for almost all sets
in the chain except $N_0$, $N_9$ and $N_{15}$. For $N_0$ we have
large size in unstable direction. This is due to the fact that $N_1$
is in some distance from $z_0$. The size $7\cdot10^{-5}$ is
necessary to reach the set $N_1$ from $N_0$.

Similar reasoning apply to $N_{15}$. The stable size must be large
enough so that the image of $N_{14}$ is captured in this direction.

On the set $N_9$ we have change of dynamics. The unstable size of
$N_9$ is smaller since on $N_8$ the parameter plays a role of
expanding direction which covers unstable direction on $N_9$. Since
this expansion is weak, we must decrease the size of $N_9$. On the
other hand the image of $N_8$ vary with the parameter and it is
larger in the stable direction of $N_9$. This forces to enlarge a
little bit the set $N_9$ in the stable direction.

Now we can state the first numerical lemma.
\begin{lemma}~\label{lem:heteroChainHenon}
The following covering relations hold true
\begin{equation*}
N_0\cover{PH}N_1\cover{PH}\cdots\cover{PH}N_{15}.
\end{equation*}
\end{lemma}
\textbf{Proof:} The assertion of the above lemma has been verified
using the interval arithmetics \cite{M} and the algorithms for
verifying the existence of covering relations described in
\cite{WZ}. We were able to verify the necessary inequalities on each
wall of $N_i$, $i=0,\ldots,15$ without any subdivision of the sets,
so the computational time was less than one second. \qed

\subsection{The cone conditions along the heteroclinic chain}

In this section we will show that the cone conditions are satisfied
for the sequence of covering relations from
Lemma~\ref{lem:heteroChainHenon}. We have two results which give us
a numerical method for verifying the cone conditions.
\begin{lemma}{\cite[Lemma~6]{KWZ}}\label{lem:coneVerification}
Let $(N,Q_N)$ and $(M,Q_M)$ be h-sets with cones in $\mathbb R^n$
and let $f:N\to\mathbb R^n$ be $\mathcal C^1$ such that
$N\cover{f}M$. Let $[Df(N)]_I$ denote the interval enclosure of
the set of matrices $Df(N)$. If the interval matrix
\begin{equation*}
V = [Df(N)]_I^T Q_M [Df(N)]_I - Q_N
\end{equation*}
is positive definite then the cone conditions are satisfied for the
covering relation $N\cover{f}M$.
\end{lemma}

Let $A = A_c + [-1,1]\cdot A_0$ be an interval matrix, where
$A_c,A_0\in\mathbb R^{n\times n}$ are real and symmetric. For
$z\in\mathbb R^n$ by $\Delta(z)$ we will denote a diagonal matrix
with $z_i$'s at the diagonal.
\begin{lemma}{\cite[Theorem~2]{R}}\label{lem:posDefVerification}
The  interval matrix $A=A_c + [-1,1]\cdot A_0$ is positive
definite if and only if for each sequence $z\in\{-1,1\}^n$ the
matrix
\begin{equation*}
A_z = A_c - \Delta(z)A_0\Delta(z)
\end{equation*}
is positive definite.
\end{lemma}
In the light of the above lemma to verify that a symmetric interval
matrix is positive definite it is enough to verify if $2^{n-1}$ real
matrices $A_z$ are positive definite.

In order to verify the cone conditions we have to define quadratic
forms on the sets $N_i$, $i=0,\ldots,15$.  Denote by
\begin{equation}\label{eq:HenonEigenvalues}
\lambda=3.858169402,\quad \mu=0.07775708341
\end{equation}
an approximate eigenvalues of $DH_{a_0}$ at $z_0$. Recall, by
$\Delta(p_1,\ldots,p_n)$ we denote a diagonal matrix with $p_i$'s at
the diagonal. For $i=0,\ldots,15$ we define the quadratic form on
the h-set $N_i$ by
\begin{equation*}
    Q_i = \Delta\left((p_i)_1,(p_i)_2,(p_i)_3,(p_i)_4\right)
\end{equation*}
where the coefficients are listed in
Table~\ref{tab:henonQuadForm}. The quadratic forms $Q_i$ are
defined in the coordinate systems given by matrixes $M_i$
(\ref{eq:HenonMatrices}) used to define the h-sets $N_i$,
$i=0,\ldots,15$. These matrices have normalized columns. In these
coordinates the sets are given by
$N_i=[(d_i)_1,\ldots,(d_i)_4]\cdot[-1,1]$, $i=0,\ldots,15$ -- see
(\ref{eq:defHSetsHenon}) and Table~\ref{tab:HenonHSets}.
\begin{table}
\begin{center}
\begin{tabular}{|c|c|c|c|c|}
    \hline
        $i$ & $(p_i)_1$& $(p_i)_2$ & $(p_i)_3$ & $(p_i)_4$\\
            & unstable dir. & stable dir. & tangent dir. & parameter\\
    \hline
        $0$ & $3/\lambda^2$ & $-\mu^2$ & $-(\mu/\lambda)^2$ & $2(1.5)^{-8}$\\
    \hline
        $1$ & $1/\lambda^2$ & $-0.1$ & $-0.5$ & $2(1.5)^{-7}$\\
    \hline
        $2$ & $1/\lambda^2$ & $-0.1$ & $-1$ & $2(1.5)^{-6}$\\
    \hline
        $3$ & $1/\lambda^2$ & $-0.1$ & $-1$ & $2(1.5)^{-5}$\\
    \hline
        $4$ & $1/\lambda^2$ & $-0.1$ & $-1$ & $2(1.5)^{-4}$\\
    \hline
        $5$ & $1/\lambda^2$ & $-0.1$ & $-1$ & $2(1.5)^{-3}$\\
    \hline
        $6$ & $1/\lambda^2$ & $-0.1$ & $-1$ & $2(1.5)^{-2}$\\
    \hline
        $7$ & $1/\lambda^2$ & $-0.1$ & $-1$ & $2(1.5)^{-1}$\\
    \hline
        $8$ & $0.5/\lambda^2$ & $-1$ & $-1$ & $2$\\
    \hline
        $9$ & $100/\lambda^2$ & $-0.1$ & $100(\mu/\lambda)^2$ & $-2$\\
    \hline
        $10$ & $40/\lambda^2$ & $-0.1$ & $(\mu/\lambda)^2$ & $-2(1.5)^{-1}$\\
    \hline
        $11$ & $10/\lambda^2$ & $-0.1$ & $(\mu/\lambda)^2$ & $-2(1.5)^{-2}$\\
    \hline
        $12$ & $1/\lambda^2$ & $-0.1$ & $(\mu/\lambda)^2$ & $-2(1.5)^{-3}$\\
    \hline
        $13$ & $1/\lambda^2$ & $-0.1$ & $(\mu/\lambda)^2$ & $-2(1.5)^{-4}$\\
    \hline
        $14$ & $1/\lambda^2$ & $-0.1$ & $(\mu/\lambda)^2$ & $-2(1.5)^{-5}$\\
    \hline
        $15$ & $0.3/\lambda^2$ & $-0.1$ & $(\mu/\lambda)^2$ & $-2(1.5)^{-6}$\\
    \hline
\end{tabular}
\end{center}
\caption{Coefficients of the quadratic forms on the sets
$N_0,\ldots,N_{15}$, where $\lambda=3.858169402$,
$\mu=0.07775708341$ are approximate eigenvalues of
$DH_{a_0}(z_0)$.\label{tab:henonQuadForm}}
\end{table}
Let us comment briefly about the choice of these coefficients.
Assume that we would like to verify the cone conditions for the
covering relation $N\cover{f}M$. Assume that $f$ is linear and in
some coordinate systems on $N$ and $M$ it is given by
$f=\Delta(\lambda_1,\ldots,\lambda_k)$. In general case, we usually
have $Df(N)$ close to a diagonal matrix, but the arguments apply.
Assume also that the quadratic forms on both sets $N$ and $M$ are
diagonal and given by
\begin{eqnarray*}
Q_N&=&\Delta(\alpha_1^N,\ldots,\alpha_k^N),\\
Q_M&=&\Delta(\alpha_1^M,\ldots,\alpha_k^M).
\end{eqnarray*}
According to Lemma~\ref{lem:coneVerification}, the cone conditions
will be satisfied if the interval matrix
\begin{equation*}
V=[Df(N)]_I^T\cdot Q_M\cdot [Df(N)]_I - Q_N =
\Delta(\lambda_1^2\alpha_1^M-\alpha_1^N,\ldots,\lambda_k^2\alpha_k^M-\alpha_k^N)
\end{equation*}
is positive definite. We see, that if $|\alpha_i^{N,M}|=1$ and
some $|\lambda_i|\gg 1$ then the corresponding coefficient in the
matrix $V$ becomes very large while for $|\lambda_j|\ll 1$ the
corresponding coefficient in $V$ is close to $1$. In general
(nonlinear) case we have some nonzero intervals off the diagonal
of $V$. Therefore, from the computational point of view, it is
better to make the matrix $V$ somehow uniform, i.e. such that the
coefficients on the diagonal are of the same magnitude. This can
be achieved by setting coefficients $\alpha_i\approx
\lambda_i^{-2}$ for these $i$ such that $|\lambda_i|\gg 1$.

This can be seen in the first column of the
Table~\ref{tab:henonQuadForm} and the second part of the third
column.

The coefficients for $N_9$ are chosen to be able to switch from
the unstable to the stable manifolds.

Notice also a different choice of the coefficients for the
quadratic form $Q_0$ in $N_0$. The reason is that we have to prove
that the center-unstable manifold at the fixed point $c_0$ is a
horizontal disk in $N_0$ satisfying the cone conditions. Therefore
we need to apply the method described in
Theorem~\ref{thm:coneInvMan} to the inverse of $PH$. In this case,
the coefficients in the second and the third column correspond to
expanding directions with approximate eigenvalues $1/\mu$ and
$\lambda/\mu$, both greater than $1$.

We have
\begin{lemma}\label{lem:coneChainHenon}
The cone conditions are satisfied for the sequence of covering
relations
\begin{equation*}
N_0\cover{PH}N_1\cover{PH}\cdots\cover{PH}N_{15}.
\end{equation*}
\end{lemma}
\textbf{Proof:} Observe that verification of the cone conditions
for these covering relations require computation of the $DPH$
which involves the second order derivatives of $H$. In the
computer assisted proof we computed an enclosure for $DPH(N_i)$,
$i=0,1\ldots,14$ using whole $N_i$ without subdivision as an
initial condition of the routine which computes $DPH$. Then we
applied Lemma~\ref{lem:coneVerification} and
Lemma~\ref{lem:posDefVerification} to prove that the cone
conditions are satisfied. The C++ program which verifies the
assertion executes within less than one second on a laptop-type
computer.\qed

\subsection{Parameterization of center-unstable and center-stable
manifolds at $c_0$ and $c_{15}$, respectively.}
\label{subsec:centerman}

 In this
section we will use the method described in
Theorem~\ref{thm:coneInvMan} in order to parameterize the
center-unstable and center-stable manifolds as a horizontal and
vertical discs satisfying the cone conditions in $N_0$ and
$N_{15}$, respectively. This, together with
Lemma~\ref{lem:heteroChainHenon} and
Lemma~\ref{lem:coneChainHenon} will give us a proof of
Theorem~\ref{thm:mainHenon}.

Put
\begin{eqnarray*}
\widetilde{N_{0}} &=& (z_0,[u_0]) + \widetilde{
M}\cdot\left((d_{0})_1,(d_{0})_2,(d_{0})_3\right),\\
\widetilde{N_{15}} &=& (z_0,[s_0]) + \widetilde{
M}\cdot\left((d_{15})_1,(d_{15})_2,(d_{15})_3\right),
\end{eqnarray*}
where $\widetilde{M}$ is a $3\times3$ minor of the matrix
$M_{15}=M_0$ (see (\ref{eq:HenonMatrices})) after removing the last
column and the last row. The coefficients $(d_i)_j$ are listed in
the Table~\ref{tab:HenonHSets}. Geometrically, these sets are just
projection onto $(x,y,t)$ coordinates of $N_0$ and $N_{15}$,
respectively.

The set $\widetilde{N_{15}}$ is a three-dimensional h-set with two
expanding directions (corresponding to eigenvalues $\lambda$ and
$\frac{\lambda}{\mu}$) and one nominally stable direction
(corresponding to the eigenvalue $\mu$), where $\lambda,\mu$ are
eigenvalues of $DH_{a_0}(z_0)$ -- see (\ref{eq:HenonEigenvalues}).

On the set $\widetilde{N_0}$ we will compute the inverse map of
$PH_a$ so that the role of nominally stable and nominally unstable
directions interchange. Hence, the set $\widetilde{N_0}$ has two
nominally unstable directions (with eigenvalues $\mu^{-1}$ and
$\lambda/\mu$) and one nominally stable direction (with eigenvalue
$\lambda^{-1}$).

On both sets we set quadratic forms defining the cones to be equal
to
\begin{eqnarray}
\widetilde{Q_{0}} &=&
-\Delta\left((p_{0})_1,(p_{0})_2,(p_{0})_3\right)\label{eq:tildeQ0},\\
\widetilde{Q_{15}} &=&
\Delta\left((p_{15})_1,(p_{15})_2,(p_{15})_3\right)\label{eq:tildeQ15},
\end{eqnarray}
where the coefficients are listed in the
Table~\ref{tab:henonQuadForm}.

\begin{lemma}\label{lem:centerStableHenon}
The center-stable manifold of $PH$ at $c_{15}$ can be parameterized
as a vertical disk in $N_{15}$ satisfying the cone conditions with
respect to the quadratic form $Q_{15}$.
\end{lemma}
\textbf{Proof:} We use Theorem~\ref{thm:coneInvMan}. Let $\pi_a$
denote a projection onto the parameter coordinate. With a computer
assistance we verified that for $a\in\pi_a(N_{15})$ holds
\begin{equation*}
\widetilde{N_{15}}\cover{PH_a}\widetilde{N_{15}}
\end{equation*}
and the cone conditions are satisfied with the quadratic form
$\widetilde{Q_{15}}$. Hence, the first assumption of
Theorem~\ref{thm:coneInvMan} is fulfilled.

Then we computed the constants $A$, $M$, $L$ which appear in
assumptions 2-3 of Theorem~\ref{thm:coneInvMan} and we obtained
\begin{eqnarray*}
A &\geq& 0.099394300936541294,\\
M & \leq & 0.084042214456891598,\\
L & \leq & 0.0070394636406844067.
\end{eqnarray*}
Hence constant $\Gamma$ from assumption 4 of
Theorem~\ref{thm:coneInvMan} can be chosen to be equal to
$\Gamma=0.57737423322563175$. With this $\Gamma$ the coefficient
$\delta$ defined in (\ref{eq:defDeltaCoeff}) is equal to
\begin{equation*}
\delta=\frac{\Gamma^2}{\|\alpha\|} \leq 16.54078540195168.
\end{equation*}
where $\alpha$ appears in the decomposition of the quadratic form
$\widetilde{Q_{15}}(x,y,t)=\alpha(x,t)-\beta(y)$ and
$\|\alpha\|=\max\left\{0.3/\lambda^2,(\mu/\lambda)^2\right\}=0.3/\lambda^2$
-- see Table~\ref{tab:henonQuadForm}.

 From Theorem~\ref{thm:coneInvMan} it follows
that the center-stable manifold of $PH$ at $c_{15}$ can be
parameterized as a vertical disk in $N_{15}$ satisfying the cone
conditions with respect to the quadratic form
\begin{equation*}
\overline{Q_{15}}(x,y,t,a) = \delta \widetilde{Q_{15}}(x,y,z) - a^2.
\end{equation*}
Recall  that by $Q_{15}$ we denote the quadratic form on $N_{15}$.
To finish the proof let us observe that
\begin{equation*}
Q_{15}(x,y,t,a) = \widetilde{Q_{15}}(x,y,t) - 2(1.5)^{-6}a^2
\end{equation*}
see (\ref{eq:tildeQ15}) and Table~\ref{tab:henonQuadForm}. Moreover,
we have $2(1.5)^{-6}\delta>1$. This shows that
\begin{equation*}
\delta Q_{15}(x,y,t,a) = \delta\left(\widetilde{Q_{15}}(x,y,t) -
2(1.5)^{-6}a^2\right) < \overline{Q_{15}}(x,y,t,a).
\end{equation*}
Therefore the center-stable manifold of $PH$ at $c_{15}$ is a
vertical disk in $N_{15}$ satisfying the cone condition for the
quadratic form $Q_{15}$.
\qed

We have a similar lemma about the parameterization of the
center-unstable manifold of $PH$ at $c_0$ as a horizontal disk in
$N_0$.
\begin{lemma}\label{lem:centerUnstableHenon}
The center-unstable manifold of $PH$ at $c_0$ can be parameterized
as a horizontal disk in $N_0$ satisfying the cone conditions with
respect to the quadratic form $Q_0$.
\end{lemma}
\textbf{Proof:} We will proceed as in
Lemma~\ref{lem:centerStableHenon} but for the map $PH^{-1}$. With a
computer assistance we verified that $a\in\pi_a(N_0)$ holds
\begin{equation*}
\widetilde{N_0}\cover{PH_a^{-1}}\widetilde{N_0}
\end{equation*}
and the cone conditions are satisfied. We computed the constants
which appear in (\ref{eq:GammaQn}--\ref{eq:defDeltaCoeff}) and we
got
\begin{eqnarray*}
A &\geq& 0.1877584261322994,\\
M & \leq & 0.2795983187542756,\\
L & \leq & 0.015049353557694945,\\
\Gamma &=&0.33278415598142302,\\
\delta & \leq & 18.316620936531205.
\end{eqnarray*}
Hence, the center-stable manifold for $PH^{-1}$ at $c_0$ is a
vertical disk in $N_0$ satisfying the cone condition for the
quadratic form
\begin{equation*}
\overline{Q_0(x,y,t,a)} = \delta\widetilde{Q_0}(x,y,t) - a^2.
\end{equation*}
From (\ref{eq:tildeQ0}) and Table~\ref{tab:henonQuadForm} we have
\begin{equation*}
Q_0(x,y,t,a) = -\widetilde{Q_0}(x,y,t,a) + 2(1.5)^{-8}a^2.
\end{equation*}
Since $2(1.5)^{-8}\delta>1$ we have
\begin{equation*}
-\delta Q_0(x,y,t,a) = \delta\left(\widetilde{Q_0}(x,y,t,a) -
2(1.5)^{-8}a^2\right) < \overline{Q_0(x,y,t,a)}
\end{equation*}
and the center-stable manifold of $PH^{-1}$ at $c_0$ can be
parameterized as a vertical disk in $N_0$ satisfying the cone
conditions with respect to the quadratic form $-Q_0$. This means
that the center-unstable manifold of $PH$ can be parameterized as a
horizontal disk in $N_0$ satisfying the cone conditions with respect
to the quadratic form $Q_0$.
\qed

\noindent\textbf{Proof of Theorem~\ref{thm:mainHenon}:} We obtain
our conclusion from Lemmas~\ref{lem:coneChainHenon},
\ref{lem:centerStableHenon} and \ref{lem:centerUnstableHenon}
combined with the discussion of the strategy of the proof in the
first part of Section~\ref{sec:covrel}. \qed

\section{Application to the forced damped pendulum equation.}
\label{sec:pendulum}

Let us consider an equation for the forced damped pendulum motion
\begin{equation}\label{eq:pendulum}
    \ddot{x} + \beta\dot x + \sin(x) = \cos(t).
\end{equation}
The equation (\ref{eq:pendulum}) is a well known example of an
equation which exhibits chaotic dynamics see \cite{GH} and
references given there. In \cite{BCGH} it has been proven (a
computer assisted proof) that for the parameter $\beta=0.1$ the
$2\pi$-time map is semiconjugated to the full shift on three
symbols on some compact invariant set.

The equation (\ref{eq:pendulum}) defines a flow on $\mathbb
R^2\times S^1$, where $S^1$ is a unit circle. Let us define the
Poincar\'e map $T_\beta\colon \mathbb R^2\to\mathbb R^2$ by
\begin{equation}
    T_\beta(x,\dot x) = (x(2\pi),\dot x(2\pi),
\end{equation}
where $x(t)$ is a solution of (\ref{eq:pendulum}) with the parameter
$\beta$.

The aim of this section is to prove the following theorem.
\begin{theorem}\label{thm:mainPendulum}
For all parameter values
\begin{equation*}
\beta\in\mathcal B=0.247133729485+[-1,1]\cdot 1.2\cdot 10^{-10}
\end{equation*}
there exists a hyperbolic fixed point for $T_\beta$ corresponding to
a $2\pi$ periodic solution of (\ref{eq:pendulum}). Moreover, there
exists a parameter value $\beta\in\mathcal B$ such that the map
$T_{\beta}$ has a quadratic homoclinic tangency unfolding
generically for that fixed point.
\end{theorem}

\begin{figure}[htbp]
\centerline{
    \includegraphics[height=2.2in]{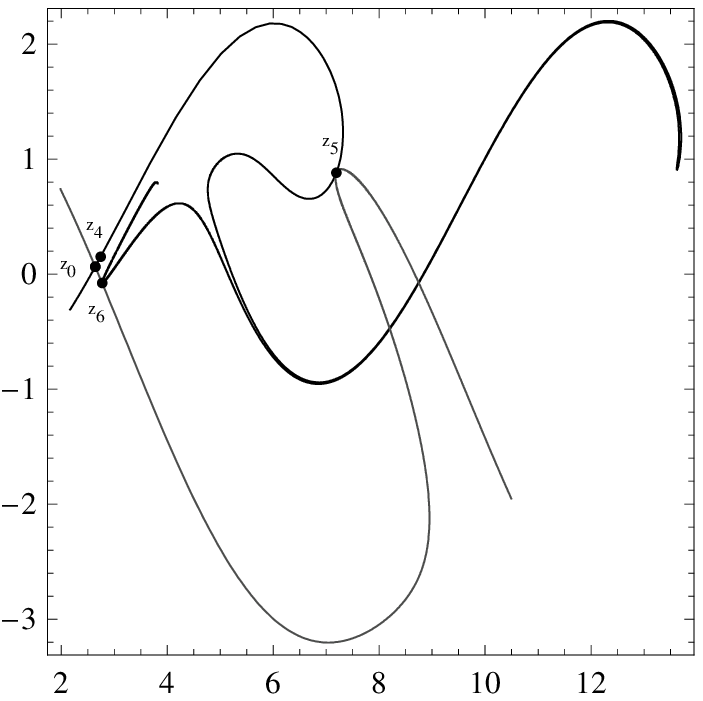}
    \includegraphics[height=2.2in]{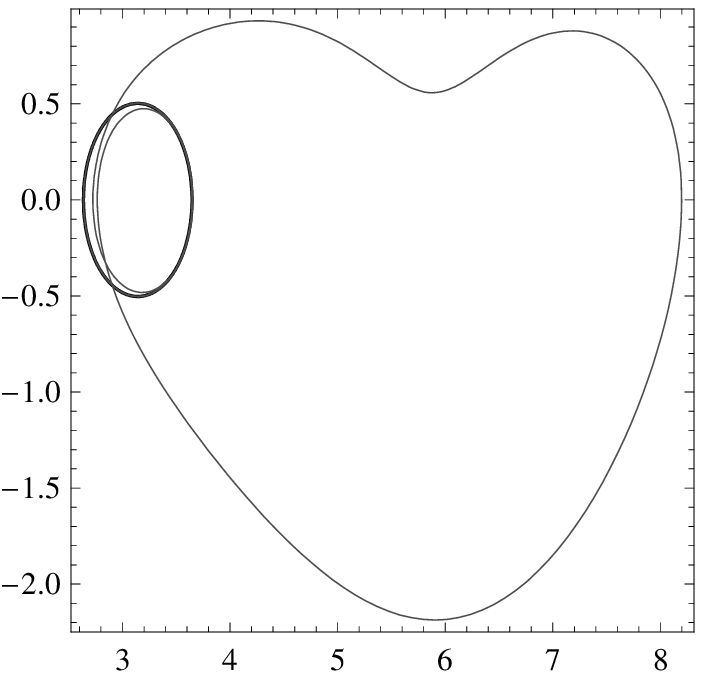}
}\caption{Left: parts of the unstable and stable manifolds of a
fixed point of $T_\beta$, for $\beta=0.247133729485$. Right:
Periodic orbit and a homoclinic orbit to the periodic orbit for
(\ref{eq:pendulum}) projected onto $(x,\dot x)$ coordinates. This
homoclinic orbit corresponds to the quadratic tangency of
invariant manifolds. The left picture suggests also the existence
of transversal homoclinic points.\label{fig:pendulumManifolds}}
\end{figure}

The proof of Theorem~\ref{thm:mainPendulum} uses the same method as
for the H\'enon map case. Since all the details has been discussed
in Section~\ref{sec:henon} we will give here only the definition of
the sets, quadratic forms and we will state the numerical lemmas.

The main difference between the H\'enon map case and that of the
forced damped pendulum equation is how we compute derivatives of a
map up to the second order, required in the proof. Since the H\'enon
map is given explicitly these derivatives can be computed by hand or
using the Automatic Differentiation tools \cite{G}. For the forced
damped pendulum equation we used the $C^r$-Lohner algorithm
presented in \cite{WZ2}. This algorithm allows us to integrate
efficiently the variational equations for ODE's and compute partial
derivatives of Poincar\'e maps.

\subsection{Heteroclinic chain of covering relations.}
Let $PT_\beta : P\mathbb R^2\to P\mathbb R^2$ denote the map induced
by $T_\beta$ on the projective bundle and let $\psi$ be a local
parameterization of the space $P\mathbb R^2$ as defined in
(\ref{eq:atlas}). Let
\begin{eqnarray*}
    PT : P\mathbb R^2\times \mathbb R\ni (x,\dot x, t,\beta) = \left((\psi^{-1}\circ PT_{\beta}\circ
    \psi)(x,\dot x,t),\beta\right)\in P\mathbb R^2\times \mathbb R.
\end{eqnarray*}

Let
\begin{eqnarray}
    \beta_0 &=& 0.247133729485,\nonumber\\
    z_0 &=& (2.6410109874338904,0.063471204982120187),\nonumber\\
    u_0 &=& (0.76818278871270265,0.64023058590290362),\label{eq:penUnstable}\\
    s_0 &=& (0.67655372773981033,-0.73639327365298942),\label{eq:penStable}\\
    \lambda &=& 211.83022271012155,\label{eq:pendulumLambda}\\
    \mu &=& 0.00099918347695168025\label{eq:pendulumMu}.
\end{eqnarray}
From numerical experiments we got that $z_0$ (see
Fig.~\ref{fig:pendulumManifolds}) is an approximate fixed point for
$T_{\beta_0}$ with approximate eigenvalues $\lambda$, $\mu$ and
approximate eigenvectors $u_0$, $s_0$, respectively. This point has
been found by the standard Newton method.

Let us denote
\begin{equation*}
    \begin{array}{rclrcl}
    &&&c_0&=&\psi^{-1}((z_0,u_0),\beta_0),\\
    z_1 &=& z_0+1.41442890240556\cdot 10^{-8}u_0, & c_1&=&\psi^{-1}((z_1,u_0),\beta_0)\\
    z_{i} &=& T_{\beta_0}(z_{i-1}),\text{ for }i=2,3,4,5, & c_{i}&=&PT(c_{i-1}),\text{ for }i=2,3,4,5,\\
    z_8 &=& z_0 + 19.0992395815\cdot 10^{-8}s_0,& c_8&=&\psi^{-1}((z_8,s_0),\beta_0)\\
    z_i &=& T_{\beta_0}^{-1}(z_{i+1})\text{ for }i=7,6, & c_{i}&=&PT^{-1}(c_{i-1})\text{ for }i=7,6,\\
    z_9 &=& z_0, & c_9 &=& \psi^{-1}((z_0,s_0),\beta_0).
    \end{array}
\end{equation*}
Some of these points are shown in Fig.~\ref{fig:pendulumManifolds}.
The points $z_1$ and $z_8$ are chosen close to a heteroclinic
trajectory for $T_ {\beta_0}$. Numerical simulation shows that
\begin{equation}\label{eq:accuracyPendulum}
PT^{-1}(c_6)-c_5\approx (3.42\cdot 10^{-12},6.34\cdot
10^{-12},-2.32\cdot 10^{-8},0).
\end{equation}
In fact, we observe that between $c_5$ and $c_6$ the role of the
tangent direction changes and at $c_6$ it becomes repelling.

The reason for which we define two points $z_1$ and $z_8$ and
compute their forward and backward trajectory, respectively, is
due to the numerical problems with forward propagation of tangent
coordinate at $z_6$, $z_7$, even in simulation only. Recall, the
tangent coordinate $t$ becomes strongly repelling  at $z_6$ and
$z_7$ with eigenvalue of the order $10^5$. Therefore, it is easier
to propagate the preimage of the stable direction at $z_8$ which
is attracting at $z_7$ and $z_6$ for the inverse map.

Approximate estimation (\ref{eq:accuracyPendulum}) shows us that the
points $c_i$, $i=1,\ldots,8$ are close to the possible existing
heteroclinic trajectory for $PT$.

Now we will define the coordinate systems of the h-sets centered at
$c_i$'s. Let $M_i$ be the matrix of coordinate system of the set
centered at $c_i$. We assume $M_i$ has the form
(\ref{eq:HenonMatrices}). These matrices are computed as follows
\begin{itemize}
\item put $u_9=u_0,\ u_1=u_0,\ u_8=u_0$ and $s_9=s_0,\ s_1=s_0,\ s_8=s_0$,
where $u_0,s_0$ are defined by
(\ref{eq:penUnstable}--\ref{eq:penStable}),
\item for $i=2,3,4,5$ we set
\begin{eqnarray*}
    u_i&=&\pi_t(c_i),\\
    s_i&=&\pi_t(PT_{\beta_0}^{-1}(z_{i+1},\pi_t(c_{i+1})^\perp)), i=2,3,4,\\
    s_5&=&u_5^\perp,
\end{eqnarray*}
\item for $i=6,7$ we set
\begin{eqnarray*}
    s_i&=&\pi_t(c_i),\\
    u_i&=&\pi_t\left(PT_{\beta_0}\left(T^{-1}(z_i),\pi_t\left(PT_{\beta_0}^{-1}(z_i,s_i)\right)^\perp\right)\right).
\end{eqnarray*}
\end{itemize}

As in the case of the H\'enon map the h-sets $N_0,\ldots,N_9$ for
the map $PT$ will be defined by formula (\ref{eq:defHSetsHenon})
with matrices $M_i$ as above and the diameters given in
Table~\ref{tab:pendulumHSets}.

\begin{table}
\begin{center}
\begin{tabular}{|c|c|c|c|c|}
    \hline
        $i$ & $10^{10}\cdot (d_i)_1$& $10^{10}\cdot (d_i)_2$ & $10^{8}\cdot (d_i)_3$ & $10^{10}\cdot (d_i)_4$\\
            & unstable dir. & stable dir. & tangent dir. & parameter\\
    \hline
        $0$ & $5$ & $0.2$ & $0.26$ & $1.2(1.01)^5$\\
    \hline
        $1$ & $0.4$ & $6$ & $2$ & $1.2(1.01)^4$\\
    \hline
        $2$ & $0.4$ & $6$ & $4$ & $1.2(1.01)^3$\\
    \hline
        $3$ & $0.6$ & $6$ & $2$ & $1.2(1.01)^2$\\
    \hline
        $4$ & $5$ & $2$ & $1$ & $1.2(1.01)$\\
    \hline
        $5$ & $27$ & $2$ & $1$ & $1.2$\\
    \hline
        $6$ & $0.4$ & $15$ & $0.4$ & $1.2(1.01)$\\
    \hline
        $7$ & $0.4$ & $8$ & $0.4$ & $1.2(1.01)^2$\\
    \hline
        $8$ & $0.4$ & $10$ & $1.2$ & $1.2(1.01)^3$\\
    \hline
        $9$ & $0.4$ & $5$ & $0.2$ & $1.2(1.01)^4$\\
    \hline
\end{tabular}
\end{center}
\caption{Diameters of the h-sets in the heteroclinic sequence for
the $PT$ map.  The h-sets $N_0,\ldots,N_5$ have two unstable
directions given by $u_i$ and the parameter. The h-sets
$N_6,\ldots,N_9$ have two unstable directions given by $u_i$ and the
tangent coordinate. \label{tab:pendulumHSets}}
\end{table}

The h-sets $N_0,\ldots,N_5$ have two unstable directions given by
$u_i$ and the parameter. The h-sets $N_6,\ldots,N_9$ have two
unstable directions given by $u_i$ and the tangent coordinate.

Now, we can state the following numerical lemma
\begin{lemma}
    The map $PT$ is well defined and continuous on $\bigcup_{i=0}^9
    N_i$. Moreover, the following covering relations hold true
\begin{equation}\label{eq:pendulumCovChain}
N_0\cover{PT}N_1\cover{PT}\ldots\cover{PT}N_9.
\end{equation}
\end{lemma}
\textbf{Proof:} In the computer assisted proof of the above lemma we
used the $C^1$-Lohner algorithm \cite{Z} and the CAPD library
\cite{CAPD} in order to integrate the variational equations for
(\ref{eq:pendulum}) and to compute the map $PT$. We used the Taylor
method of the order $20$ with a variable time step. To verify the
inequalities required for the covering relations we
 subdivided  the boundary
of each h-set with a grid depending on the set under consideration.
The total number of boxes we used is $5546$. The C++ program which
verifies the existence of covering relations
(\ref{eq:pendulumCovChain}) executes within 18 seconds on a computer
with the Intel Xeon 5160, 3GHz processor.
\qed

\subsection{The cone conditions along the heteroclinic chain of covering relations for the map $PT$.}
Recall, that for a $p\in\mathbb R^n$ by $\Delta(p)$ we denoted a
diagonal matrix with $p_i$'s on the diagonal.

For $i=0,\ldots,9$ we define the quadratic form on the h-set $N_i$
by
\begin{equation*}
    Q_i = \Delta\left((p_i)_1,(p_i)_2,(p_i)_3,(p_i)_4\right)
\end{equation*}
where the coefficients are listed in
Table~\ref{tab:pendulumQuadForm}.

Let us comment about the choice of these coefficients. In the
example presented in Section~\ref{sec:example} there is no
dependency on the parameter, hence the cone conditions are easily
achievable. In general case this dependency has huge influence for
the choice of the parameters in quadratic forms corresponding to the
parameter variable. In fact, this sometimes forces large scaling
(not only by some small factor), like in the H\'enon map case).

The other constraints on the coefficients are related to the
parameterization of center-unstable and center-stable manifolds in
the first and last sets in the heteroclinic chain of covering
relations.

The main constraint, however, appears for the covering relation in
the switch between manifolds. Here it is necessary to set
relatively large coefficients corresponding to both unstable
variables in the set after this switch and small coefficient for
the stable variable in the main phase space. This can be seen in
Tables~\ref{tab:henonQuadForm} and \ref{tab:pendulumQuadForm} for
the sets $N_9$ and $N_6$, respectively. In the next sets we can
use the hyperobolicity of the map to make these coefficients more
uniform and to reach constraints at the begin and at the end of
the chain.

\begin{table}
\begin{center}
\begin{tabular}{|c|c|c|c|c|}
    \hline
        $i$ & $(p_i)_1$& $ (p_i)_2$ & $(p_i)_3$ & $(p_i)_4$\\
            & unstable dir. & stable dir. & tangent dir. & parameter\\
    \hline
        $0$ & $80/\lambda^2$ & $-\mu^2$ & $-(\mu/\lambda)^2$ & $(1.1)^{-5}$\\
    \hline
        $1$ & $1/\lambda^2$ & $-0.01$ & $-10^{-7}$ & $(1.1)^{-4}$\\
    \hline
        $2$ & $1/\lambda^2$ & $-1$ & $-10^{-5}$ & $(1.1)^{-3}$\\
    \hline
        $3$ & $1/\lambda^2$ & $-1$ & $-10^{-5}$ & $(1.1)^{-2}$\\
    \hline
        $4$ & $1/\lambda^2$ & $-1$ & $-10^{-5}$ & $(1.1)^{-1}$\\
    \hline
        $5$ & $10/\lambda^2$ & $-1$ & $-10^{-5}$ & $3$\\
    \hline
        $6$ & $1000/\lambda^2$ & $-10^{-4}$ & $10^{6}(\mu/\lambda)^2$ & $-(1.1)^{-1}$\\
    \hline
        $7$ & $1/\lambda^2$ & $-1$ & $(\mu/\lambda)^2$ & $-(1.1)^{-2}$\\
    \hline
        $8$ & $1/\lambda^2$ & $-1$ & $(\mu/\lambda)^2$ & $-(1.1)^{-3}$\\
    \hline
        $9$ & $1/\lambda^2$ & $-1$ & $(\mu/\lambda)^2$ & $-(1.1)^{-4}$\\
    \hline
\end{tabular}
\end{center}
\caption{Coefficients of the quadratic forms on the sets
$N_0,\ldots,N_9$, where $\lambda$, $\mu$ are approximate eigenvalues
of $DT_{\beta_0}(z_0)$ and are given in
(\ref{eq:pendulumLambda}-\ref{eq:pendulumMu}).\label{tab:pendulumQuadForm}}
\end{table}

For $i=2,\ldots,9$ the coordinate systems on $N_i$ are given by
the matrices $M_i$ used to define the sets $N_i$, respectively.

We have the following numerical result
\begin{lemma}
The cone conditions are satisfied for the sequence of covering
relations
\begin{equation*}
N_0\cover{PT}N_1\cover{PT}\cdots\cover{PT}N_9.
\end{equation*}
\end{lemma}
\textbf{Proof:} By Lemma~\ref{lem:coneVerification} it is enough to
verify that for $i=0,\ldots,8$ the interval matrix
\begin{equation*}
V_i = [DPT(N_i)]_I^T Q_{N_{i+1}} [DPT(N_i)]_I - Q_{N_i}
\end{equation*}
is positive definite, where $DPT$ is computed in the coordinate
systems of $N_{i}$ and $N_{i+1}$. Notice, to compute $DPT$ we need
second order derivatives of the map $T$. We used the $C^2$-Lohner
algorithm \cite{WZ2} and the CAPD library \cite{CAPD} in order to
integrate the second order variational equations for
(\ref{eq:pendulum}). We used the Taylor method of the order $20$ and
a variable time step. No subdivision of the sets $N_i$ were
necessary, i.e. whole sets $N_i$ were used as an initial condition
for the routine which computes $DPT$. The C++ program which verifies
the cone conditions from this lemma executes within $1$ seconds a
computer with the Intel Xeon 5160, 3GHz processor. \qed

\subsection{Parameterization of center-unstable and center-stable
manifolds at $c_0$ and $c_9$, respectively.}
\begin{lemma}\label{lem:centerPendulum}
The following statements hold true.
\begin{itemize}
\item The center-stable manifold of $PT$ at $c_9$ can be parameterized
as a vertical disk in $N_9$ satisfying the cone conditions with
respect to the quadratic form $Q_9$.

\item The center-unstable manifold of $PT$ at $c_0$ can be
parameterized as a vertical disk in $N_0$ satisfying the cone
conditions with respect to the quadratic form $Q_0$.
\end{itemize}
\end{lemma}
The computer assisted proof of the above lemma is essentially the
same as the proof of Lemma~\ref{lem:centerStableHenon} and
Lemma~\ref{lem:centerUnstableHenon} for the H\'enon map. Therefore
we skip the details. The C++ program which verifies the assertion
executes within 11 seconds on the Intel Xeon 5160, 3GHz processor.
The most time-consuming part is the verification of the existence of
covering relations for the projected sets in three-dimensional
space.

We used the Taylor method of the order $20$ and a variable time
step when check the covering relations and when integrate second
order variational equations.

\vskip\baselineskip\noindent\textbf{Proof of
Theorem~\ref{thm:mainPendulum}:} We conclude the proof as in the
H\'enon map example. \qed

\section{Implementation notes.}
In order to compute bounds for the H\'enon map and the Poincar\'e
map $T_\beta$ and their derivatives we used the interval arithmetic
\cite{M}, automatic differentiation \cite{G} and the $C^r$-Lohner
algorithm \cite{WZ2} developed at the Jagiellonian University by the
CAPD group \cite{CAPD}. The C++ source files of the program with an
instruction how it should be compiled and run are available at
\cite{W}.

The program has been tested under several linux distributions,
including 32 and 64 bits architectures and the gcc compiler versions
4.1.2, 4.2.1 and 4.3.2 on the Intel Pentium IV, Intel Core 2 Duo,
Intel Xeon and the AMD Quad Core processors.

\end{document}